\documentclass[11pt]{amsart}

\usepackage{amsmath}

\usepackage{amssymb}

\newcommand{\C}{{\mathbb C}}       
\newcommand{\R}{{\mathbb R}}       
\newcommand{\DD}{{\mathcal D}}
\newcommand{\HH}{{\mathcal H}}
\newcommand{\LL}{{\mathcal L}}

\newcommand{\CE}{{{\mathcal C}_{\varepsilon}}}

\newcommand{\diam}{{\rm diam}}
\newcommand{\dist}{{\rm dist}}
\newcommand{\ds}{\displaystyle }
\newcommand{\ts}{\textstyle }

\newcommand{\lra}{{\longrightarrow}}

\newcommand{\real}{{\rm Re}}
\newcommand{\rf}[1]{{(\ref{#1})}}

\newcommand{\supp}{{\rm supp}}

\newcommand{\uuunt}{{\int\!\!\int\!\!\int}}
\newcommand{\vphi}{{\varphi}}
\newcommand{\ve}{{\varepsilon}}
\newcommand{\vv}{{\vspace{2mm}}}

\newcommand{\wt}[1]{{\widetilde{#1}}}
\newcommand{\wh}[1]{{\widehat{#1}}}

\newcommand{\bmo}{{B\!M\!O}}
\newcommand{\rbmo}{{R\!B\!M\!O}}
\newcommand{\hba}{{H^{1,\infty}_{atb}}}
\newcommand{\hbap}{{H^{1,p}_{atb}}}
\newcommand{\hbm}{{H^{1,\infty}_{atb}(\mu)}}
\newcommand{\hbp}{{H^{1,p}_{atb}(\mu)}}

\newtheorem{theorem}{Theorem}[section]
\newtheorem{lemma}[theorem]{Lemma}
\newtheorem{coro}[theorem]{Corollary}
\newtheorem{propo}[theorem]{Proposition}

\theoremstyle{definition}

\newtheorem{example}[theorem]{Example}

\theoremstyle{remark}
\newtheorem{rem}[theorem]{Remark}

\numberwithin{equation}{section}

\newcommand{\brem}{\begin{rem}}
\newcommand{\erem}{\end{rem}}

\newcommand{\bexam}{\begin{example}}
\newcommand{\eexam}{\end{example}}

\begin{document}

\title[$\bmo$ and $H^1$ for non doubling measures]
{$\bmo$, $H^1$, and Calder\'on-Zygmund operators for
non doubling measures}

\author[XAVIER TOLSA]{Xavier Tolsa}

\address{Departament de Matem\`atica Aplicada i An\`alisi, Universitat de
Barcelona, Gran Via 585, 08071 Barcelona, Spain}

\curraddr{Department of Mathematics, Chalmers, 412 96 G\"oteborg, Sweden}

\email{xavier@math.chalmers.se}

\thanks{Partially supported by grants DGICYT PB96-1183 and CIRIT
1998-SGR00052 (Spain).}


\date{February 17, 2000.}


\begin{abstract}
Given a Radon measure $\mu$ on $\R^d$, which may be non doubling, we introduce
a space of type
$\bmo$ with respect to this measure. It is shown that many properties
which hold for the classical space $\bmo(\mu)$ when $\mu$ is a doubling
measure
remain valid for the space of type $\bmo$ introduced in this paper, without
assuming $\mu$ doubling. For instance,
Calder\'on-Zygmund operators
which are bounded on $L^2(\mu)$ are also bounded from $L^\infty(\mu)$ into
the new $\bmo$ space. Moreover, this space also satisfies a John-Nirenberg
inequality, and its predual is an atomic space $H^1$. Using a sharp maximal
operator it is shown that operators which are bounded from $L^\infty(\mu)$
into the new $\bmo$ space and from its predual $H^1$ into $L^1(\mu)$ must
be bounded on $L^p(\mu)$, $1<p<\infty$. From this result one can obtain
a new proof of the $T(1)$ theorem for the Cauchy transform for non doubling
measures. Finally, it is proved that commutators of Calder\'on-Zygmund
operators bounded on $L^2(\mu)$ with functions of the new $\bmo$ are
bounded on $L^p(\mu)$, $1<p<\infty$.
\end{abstract}

\maketitle

\newpage

\section{Introduction} \label{sec1}

In this paper, given a Radon measure $\mu$ on $\R^d$ which may be non
doubling, we introduce a $\bmo$ space and an
atomic space (the predual of the $\bmo$ space), with respect to the
measure $\mu$. It is shown that, in some ways, these spaces play the role
of the classical spaces $\bmo$ and $H^1_{at}$ in case $\mu$ is a doubling
measure.

Recently it has been proved that many results of the Cader\'on-Zygmund theory
of singular integrals remain valid for non doubling measures on $\R^d$.
A version of the $T(1)$ theorem for the Cauchy transform was obtained
in \cite{Tolsa1}, and another for more general Calder\'on-Zygmund
operators in \cite{NTV1}. Cotlar's inequality and weak $(1,1)$ estimates were
studied in \cite{NTV1} (for the particular case of the Cauchy transform, the
weak $(1,1)$ estimate was studied also in \cite{Tolsa1}).
G. David \cite{David} obtained a theorem of $T(b)$ type for non doubling
measures that solved Vitushkin's conjecture for sets with
positive finite $1$-dimensional Hausdorff measure. Another $T(b)$ theorem
suitable for solving this conjecture was proved later by
Nazarov, Treil and Volberg \cite{NTV3}. In \cite{Tolsa2},
it is shown that if the Cauchy transform is bounded on $L^2(\mu)$,
then the principal values of the Cauchy transform
exist $\mu$-almost everywhere in $\C$.
In \cite{Verdera}, it is given
another proof for the $T(1)$ theorem for the Cauchy
transform, and in \cite{Tolsa3} a $T(1)$ theorem suitable
for non doubling measures {\em with atoms} is proved.
Also, in \cite{NTV*},
another $T(b)$ theorem for non doubling measures (closer to the classical one
than the ones stated above) is obtained.

However, for the moment, the attempts to find good substitutes for the space
$\bmo$ and its predual $H^1_{at}$ for non doubling measures have not been
completely succesful.
Mateu, Mattila, Nicolau and Orobitg \cite{MMNO} have studied the spaces
$\bmo(\mu)$ and $H^1_{at}(\mu)$ (with definitions similar to the classical
ones) for a non doubling measure $\mu$.
They have shown that some of the properties that these spaces satisfy when
$\mu$ is a doubling measure are satisfied also if $\mu$ is non doubling.
For example, the
John-Nirenberg inequality holds, $\bmo(\mu)$ is the dual of $H^1_{at}(\mu)$ and the
operators which are bounded from $H^1_{at}(\mu)$ into $L^1(\mu)$ and from
$L^\infty(\mu)$ into $\bmo(\mu)$ are bounded on $L^p(\mu)$, $1<p<\infty$.
Nevertheless, unlike in the case of doubling measures, Calder\'on-Zygmund
operators may be bounded on $L^2(\mu)$ but not from $L^\infty(\mu)$ into
$\bmo(\mu)$ or from $H^1_{at}(\mu)$ into $L^1(\mu)$, as it is shown by Verdera
\cite{Verdera}.
This is the main drawback of the spaces $\bmo(\mu)$ and $H^1_{at}(\mu)$
considered in \cite{MMNO}.

On the other hand, Nazarov, Treil and Volberg \cite{NTV*} have introduced
another space of $\bmo$ type. Calder\'on-Zygmund operators which
are bounded on $L^2(\mu)$ are bounded from $L^\infty(\mu)$ into their
$\bmo$ space. However, the $\bmo$ space considered in \cite{NTV*} does not
satisfy John-Nirenberg inequality, it is not known which is its predual,
and (by now) there is no any interpolation result such as the one given
in \cite{MMNO}.

\vv
Let us introduce some notation and definitions.
Let $d,n$ be some fixed integers with $1\leq n\leq d$. A kernel
$k(\cdot,\cdot)\in
L^1_{loc}(\R^d\times\R^d \setminus \{(x,y) :x=y\})$ is called a
Calder\'on-Zygmund kernel if
\begin{enumerate}
\item $\ds |k(x,y)|\leq \frac{C}{|x-y|^n},$
\item there exists $0<\delta\leq1$ such that
$$|k(x,y)-k(x',y)| + |k(y,x)-k(y,x')|
\leq C\,\frac{|x-x'|^\delta}{|x-y|^{n+\delta}}$$
if $|x-x'|\leq |x-y|/2$.
\end{enumerate}
Throughout all the paper we will assume that $\mu$ is a Radon measure on
$\R^d$ satisfying the following growth condition:
\begin{equation}  \label{creix}
\mu(B(x,r))\leq C_0\,r^n\quad \mbox{for all $x\in \R^d,\, r>0.$}
\end{equation}
The Calder\'on-Zygmund operator (CZO) associated to the kernel
$k(\cdot,\cdot)$ and the measure $\mu$ is defined (at least, formally) as
$$Tf(x) = \int k(x,y)\,f(y)\,d\mu(y).$$
The above integral may be not convergent for many functions $f$ because
the kernel $k$ may have a singularity for $x=y$. For this reason,
one introduces the truncated operators $T_\varepsilon$, $\varepsilon>0$:
$$T_\varepsilon f(x) = \int_{|x-y|>\varepsilon} k(x,y)\,f(y)\,d\mu(y),$$
and then one says that $T$ is bounded on $L^p(\mu)$ if the operators
$T_\varepsilon$ are bounded on $L^p(\mu)$ uniformly on $\varepsilon>0$.

Recall that a function $f\in L^1_{loc}(\mu)$ is said to belong to $\bmo(\mu)$
if there exists some constant $C_1$ such that
\begin{equation} \label{dbmo}
\sup_{Q} \frac{1}{\mu(Q)}\int_Q |f-m_Q(f)|\,d\mu\leq C_1,
\end{equation}
where the supremum is taken over all the cubes $Q\subset\R^d$ centered at
some point of $\supp(\mu)$ (in the paper by a cube we mean
a closed cube with sides parallel to the axes, and if $\|\mu\|<\infty$, we
allow $Q=\R^d$ too)
and $m_Q(f)$ stands for the mean of $f$ over $Q$ with respect to $\mu$, i.e.
$m_Q(f)=\int_Q f\,d\mu/\mu(Q)$. The optimal constant $C_1$ is the $\bmo$ norm
of $f$.

Let us remark that there is a slight difference between the space $\bmo(\mu)$
that we have just defined and the one considered in \cite{MMNO}: We have
taken the supremum in \rf{dbmo} over cubes which are centered at some point in
$\supp(\mu)$, while in \cite{MMNO} that supremum is taken over all the cubes
in $\R^d$.

It is well known that if $\mu$ is a doubling measure, i.e.
$\mu(2Q)\leq C\,\mu(Q)$ for all the cubes $Q$ centered at some point of
$\supp(\mu)$, and $T$ is bounded on
$L^2(\mu)$, then $T$ is also bounded from $L^\infty(\mu)$ into $\bmo(\mu)$.
As stated above, this may fail if $\mu$ is non doubling.
Hence if one wants to work with a $\bmo$ space which fulfils some of the
usual and fundamental properties related with CZO's,
then one must introduce a new space $\bmo$.
So, for a fixed $\rho>1$, one says that a function $f\in L^1_{loc}(\mu)$
belongs to $\bmo_\rho(\mu)$ if for some constant $C_2$
\begin{equation} \label{nbmo}
\sup_{Q} \frac{1}{\mu(\rho Q)}
\int_Q |f-m_Q(f)|\,d\mu\leq C_2,
\end{equation}
again with the supremum taken over all the cubes $Q$
centered at some point of $\supp(\mu)$.
This is {\em almost} the definition taken by Nazarov, Treil and Volberg in
\cite{NTV*}.

In fact, in \cite{NTV*}, the supremum in the definition \rf{nbmo} of
$\bmo_\rho$ is taken not only over all cubes $Q$ centered at some point of
$\supp(\mu)$, but over {\em all} the cubes $Q\subset\R^d$. We have prefered to
take the supremum only over cubes centered in points of $\supp(\mu)$, to be
coherent with our definitions below.

It is straigthforward to check that
if $T$ is bounded on $L^2(\mu)$, then $T$ is bounded from $L^\infty(\mu)$
into $\bmo_\rho(\mu)$. This is the  main advantage of $\bmo_\rho(\mu)$ over
$\bmo(\mu)$.
Nevertheless, the new space $\bmo_\rho(\mu)$ does not have all the nice
properties that one may expect. First of all, it happens that the
definition of $\bmo_\rho(\mu)$ depends on the constant $\rho>1$ that we
choose. Obviously, the $\bmo_\rho$ norm of $f$
(i.e. the optimal constant $C_2$ in \rf{nbmo}) depends on $\rho$. Moreover,
it is shown in \cite{NTV*} that there exist measures $\mu$ and functions
$f$ which for some $\rho>1$ are in $\bmo_\rho(\mu)$, but not for other
$\rho>1$.

Given $p\in[1,\infty)$, one says that $f\in \bmo^p_\rho(\mu)$ if
\begin{equation} \label{nbmop}
\sup_{Q} \int_Q |f-m_Q(f)|^p\,d\mu\leq C \,\mu(\rho Q).
\end{equation}
In case $\mu$ is doubling measure, by John-Nirenberg inequality, all the
spaces $\bmo^p(\mu)\equiv\bmo^p_{\rho=1}(\mu)$
coincide. This is not the case if $\mu$ is non
doubling. In \cite{NTV*} it is shown that there are measures $\mu$
and functions $f$ such that $f$ is in $\bmo^p_\rho(\mu)$ only for a proper
subset of $p\in[1,\infty)$.

In this paper we will introduce a new variant of the space $\bmo$ suitable for
non doubling measures, which will satisfy some of the properties of the
usual $\bmo$, such as for example the John-Nirenberg inequality. This
space will be a (proper, in general) subspace of the spaces
$\bmo_\rho^p(\mu)$. It will be small enough to fulfil the properties
that we have mentioned and big enough in order that CZO's
which are bounded on $L^2(\mu)$ be also bounded from
$L^\infty(\mu)$ into our new space of $\bmo$ type.

We will show that if $T$ is bounded on $L^2(\mu)$ and $g\in L^\infty(\mu)$,
then the oscillations of $f=T(g)$ satisfy not only the condition given by
\rf{nbmo}, but other regularity conditions. Then, the functions of our new
space will be the functions satisfying \rf{nbmo} and, also, these additional
regularity conditions about their oscillations. We will denote it as
$\rbmo(\mu)$ (this stands for `{\em regular} bounded mean oscillations').
Notice that we have not written $\rbmo_\rho(\mu)$. This is because, as we
will see, the definition will not depend on $\rho$, for $\rho>1$.

\vv
If one says that $f$ is in $\bmo_\rho(\mu)$ when it satisfies \rf{nbmo}, it
seems that we have to consider the atomic space $H^{1,\infty}_{at,\rho}(\mu)$
made up with functions of the form
\begin{equation} \label{hat1}
f=\sum_i \lambda_i\, a_i,
\end{equation}
where $\lambda_i\in\R$, $\sum_i|\lambda_i|<\infty$ and, for each $i$, $a_i$ is
a function supported
in a cube $Q_i$, with $\|a_i\|_{L^{\infty}(\mu)}\leq \mu(\rho Q_i)^{-1}$,
and $\int a_i\,d\mu=0$ (that is $a_i$ is an {\em atom}).
Obviously, $H^{1,\infty}_{at,\rho}(\mu)$ is the usual atomic space
$H^{1,\infty}_{at}(\mu)\equiv H^{1,\infty}_{at,\rho=1}(\mu)$
when $\mu$ is a doubling measure.
With this definition, a CZO which is bounded
in $L^2(\mu)$, is also bounded from $H^{1,\infty}_{at,\rho}(\mu)$ into
$L^1(\mu)$ (taking $\rho>1$).

In this paper we will introduce another space of atomic type:
$\hba(\mu)$ (the subindex `atb' stands for `atomic block').
This space will be made up of functions of the form
\begin{equation} \label{hat2}
f=\sum_i b_i,
\end{equation}
where the functions $b_i$ will be some elementary functions, which we will
call {\em atomic blocks} (in particular, any atom $a_i$ such as the one of
\rf{hat1} will be an atomic block).
So we will have
$H^{1,\infty}_{at,\rho}(\mu)\subset\hba(\mu)$
but, in general, $\hba(\mu)$ will be strictly bigger than
$H^{1,\infty}_{at,\rho}(\mu)$.

We will see that this new atomic space enjoys some very interesting
properties.
First of all, the definition of the space will be independent of the chosen
constant $\rho>1$. Also,
CZO's which are bounded on
$L^2(\mu)$ will be also bounded from $\hba(\mu)$ into $L^1(\mu)$.
Moreover, we will show that $\hba(\mu)$ is the predual of $\rbmo(\mu)$,
and that, as in the doubling case, there is a collection of spaces
$H^{1,p}_{atb}(\mu)$, $p>1$, that coincide with $\hba(\mu)$.

\vv
We will show two applications of all the results obtained about
$\rbmo(\mu)$ and $\hba(\mu)$. In our first application we will obtain an
interpolation theorem: We will prove
that if a linear operator is bounded from $L^\infty$ into
$\rbmo(\mu)$ and from $\hba(\mu)$ into $L^1(\mu)$, then
it is bounded on $L^p(\mu)$, $1<p<\infty$. As a consequence we will
obtain a new proof of the $T(1)$ theorem for the Cauchy transform for
non doubling measures.

We have already mentioned that in \cite{MMNO} it is also proved a theorem of
interpolation between $(H^1_{at}(\mu), L^1(\mu))$ and
$(L^\infty(\mu),\bmo(\mu))$, with $\mu$ non doubling.
However, from this result it is not possible to get the $T(1)$ theorem
for the Cauchy transform, as it is explained in \cite{MMNO}.

Finally, in our second application we will show that if a
CZO is bounded on $L^2(\mu)$, then the commutator
of this operator with a function of $\rbmo(\mu)$ is bounded on $L^p(\mu)$,
$1<p<\infty$.


\section{The space $\rbmo(\mu)$}

\subsection{Introduction}

If $\mu$ is a doubling measure and $f$ is a function belonging to $\bmo(\mu)$,
it is easily checked that if $Q,R$ are two cubes of comparable size with
$Q\subset R$, then
\begin{equation} \label{reg1}
|m_Q(f)-m_R(f)|\leq C \,\|f\|_{\bmo(\mu)}
\end{equation}
In case $\mu$ is not doubling and $f\in \bmo_\rho(\mu)$, it is easily seen
that
\begin{equation}  \label{reg2}
|m_Q(f)-m_R(f)|\leq \frac{\mu(\rho R)}{\mu(Q)} \,\|f\|_{\bmo_\rho(\mu)},
\end{equation}
and that's all one can obtain. So if $\mu(Q)$ is much smaller than $\mu(R)$,
then $m_Q(f)$ may be very different from $m_R(f)$, and one does not have any
useful information.
However, to prove most results dealing with functions in $\bmo$, some kind
of control in the changes of the mean values of $f$, such as the one in
\rf{reg1}, appears to be essential.

We will see that if $T$ is a CZO that
is bounded on $L^2(\mu)$ and $g\in L^\infty(\mu)$, then the oscillations of
$T(g)$ satisfy some properties  which will be stronger than \rf{reg2}.
Some of these properties will be stated in terms of some coefficients
$K_{Q,R}$, for $Q\subset R$ cubes in $\R^d$, which now we proceed to
describe.


\subsection{The coefficients $K_{Q,R}$}  \label{kqr}

Throughout the rest of the paper, unless otherwise stated, any
cube will be a cube in $\R^d$ {\em with sides parallel to the axes
and centered at some point of $\supp(\mu)$}.

Given two cubes $Q\subset R$ in $\R^d$, we set
\begin{equation} \label{kq1}
K_{Q,R} = 1+ \sum_{k=1}^{N_{Q,R}} \frac{\mu(2^kQ)}{l(2^kQ)^n},
\end{equation}
where $N_{Q,R}$ is the first integer $k$ such that $l(2^kQ)\geq l(R)$
(in case $R=\R^d\neq Q$, we set $N_{Q,R} = \infty$).
The coefficient $K_{Q,R}$ measures how close $Q$ is to $R$, in some sense.
For example, if $Q$ and $R$ have comparable sizes, then $K_{Q,R}$ is bounded
above by some constant which depends on the ratio $l(R)/l(Q)$ (and on the
constant $C_0$ of \rf{creix}).

Given $\alpha>1$ and $\beta>\alpha^n$, we say that some cube $Q\subset \R^d$
is $(\alpha,\beta)$-doubling if $\mu(\alpha Q) \leq \beta\,\mu(Q)$.
Due to the fact that $\mu$ satisfies the growth condition \rf{creix}, there
are a lot of ``big'' doubling cubes. To be precise, given any point
$x\in\supp(\mu)$ and $d>0$, there exists some
$(\alpha,\beta)$-doubling cube $Q$ centered at $x$ with $l(Q)\geq d$. This is
easily seen by the growth condition \rf{creix} for $\mu$ and the fact
that $\beta>\alpha^n$.

In the following lemma we show some of the properties of the
coefficients $K_{Q,R}$.

\begin{lemma} \label{kq}
We have:
\begin{enumerate}
\item If $Q\subset R\subset S$ are cubes in $\R^d$, then
$K_{Q,R}\leq K_{Q,S}$, $K_{R,S}\leq C\,K_{Q,S}$ and
$K_{Q,S}\leq C\,(K_{Q,R}+K_{R,S})$.

\item If $Q\subset R$ have comparable sizes, $K_{Q,R}\leq C$.

\item If $N$ is some positive integer and the cubes
$2Q,\,2^2Q,\ldots2^{N-1}$ are non $(2,\beta)$-doubling (with $\beta>2^n$),
then
$K_{Q,2^NQ}\leq C$, with $C$ depending on $\beta$ and $C_0$.

\item If $N$ is some positive integer and for some $\beta< 2^n$,
$$\mu(2^NQ)\leq \beta \mu(2^{N-1} Q)\leq \beta^{2} \mu(2^{N-2}Q)\leq
\ldots \leq \beta^{N}\mu(Q),$$
then $K_{Q,2^NQ}\leq C$, with $C$ depending on $\beta$ and $C_0$.
\end{enumerate}
\end{lemma}

\begin{proof} The properties 1 and 2 are immediate. Let us see 3.
For $\beta > 2^n$, we have
$\mu(2^{k+1}Q)>\beta \,\mu(2^kQ)$ for $k=1,\ldots,N-1$. Thus
$$\mu(2^k Q)< \frac{\mu(2^N Q)}{\beta^{N-k}}$$
for $k=1,\ldots,N-1$. Therefore,
\begin{eqnarray*}
K_{Q,2^NQ}& \leq& 1+ \sum_{k=1}^{N-1} \frac{\mu(2^N
Q)}{\beta^{N-k}\,l(2^kQ)^n} + \frac{\mu(2^N Q)}{l(2^NQ)^n} \\
&\leq& 1+C_0+ \frac{\mu(2^N Q)}{l(2^NQ)^n}
\sum_{k=1}^{N-1} \frac{1}{\beta^{N-k}\,2^{(k-N)n}} \\
& \leq & 1 +C_0+ C_0 \,\sum_{k=1}^\infty (2^n/\beta)^k \leq C.
\end{eqnarray*}

Let us check the fourth property. For $\beta < 2^n$, we have
\begin{eqnarray*}
K_{Q,2^NQ}& \leq& 1+ \sum_{k=1}^{N} \frac{\beta^{k}\mu(Q)}{l(2^kQ)^n}\\
&\leq & 1+ \frac{\mu(Q)}{l(Q)^n} \sum_{k=1}^{N}
\frac{\beta^{k}}{2^{kn}}\\
&\leq & 1+C_0 \sum_{k=1}^\infty \left(\frac{\beta}{2^n}\right)^k \leq C.
\end{eqnarray*}
\end{proof}

Notice that, in some sense, the property 3 of Lemma \ref{kq} says that
if the density of the measure $\mu$ in the concentric cubes grows
much faster than the size of cubes, then the coefficients
$K_{Q,2^NQ}$ remain bounded, while the fourth property says that if the
measure grows too slowly, then they also remain bounded.

\brem  \label{rem1}
If we substitute the numbers $2^k$ in the definition \rf{kq1} by
$\alpha^k$, for some $\alpha>1$, we will obtain a coefficient
$K^\alpha_{Q,R}$. It is easy to check that $K_{Q,R}\approx K^\alpha_{Q,R}$
(with constants that may depend on $\alpha$ and $C_0$).

Also, if we set
$$K_{Q,R}'= 1 + \int_{l(Q)}^{l(R)} \frac{\mu(B(x_Q,r))}{r^{n-1}}\,dr,$$
or
$$K_{Q,R}''= 1 + \int_{l(Q)\leq |y-x_Q| \leq l(R)} 
\frac{1}{|y-x_Q|^{n}}\,d\mu(y),$$
where $x_Q$ is the center of $Q$, then it is easily seen that
$K_{Q,R}\approx K_{Q,R}'\approx K_{Q,R}''$. The definitions of $K_{Q,R}'$ and
$K_{Q,R}''$ have the advantage
of not depending on the grid of cubes, unlike the one of $K_{Q,R}$.
\erem

We have stated above that there a lot of ``big'' $(\alpha,\beta)$-doubling
cubes.
In the next remark we show that, for $\beta$ big enough, there are also many
``small'' $(\alpha,\beta)$-doubling cubes.

\brem \label{rem9}
Given $\alpha>1$, if $\mu$ is any Radon measure on $\R^d$, it is known that
for $\beta$ big enough (depending on $\alpha$ and $d$), for
$\mu$-almost all
$x\in\R^d$ there is a sequence of {\em $(\alpha,\beta)$-doubling} cubes
$\{Q_n\}_n$ centered at $x$ with $l(Q_n)$ tending to 0 as $n\to\infty$.

For $\alpha=2$, we denote by $\beta_d$ one of these big constants $\beta$.
For definiteness, one can assume that $\beta_d$ is twice the infimum of
these $\beta$'s.

If $\alpha$ and $\beta$ are not specified, by a doubling cube we will
mean a $(2,\beta_d)$-doubling cube.

Let $f\in L^1_{loc}(\mu)$ be given.
Observe that, by the Lebesgue differentiation theorem, for $\mu$-almost all
$x\in\R^d$ one can find
a sequence of $(2,\beta_d)$-doubling cubes $\{Q_k\}_k$
centered at $x$ with $l(Q_k)\to 0$ such that
$$\lim_{k\to\infty} \frac{1}{\mu(Q_k)} \int_{Q_k} f \,d\mu = f(x).$$
Thus, for any fixed $\lambda>0$, for $\mu$-almost all $x\in\R^d$ 
such that
$|f(x)|>\lambda$, there exists a sequence of $(2,\beta_d)$-doubling 
cubes $\{Q_k\}_k$ centered at $x$ with $l(Q_k)\to 0$ such that
$$\limsup_{k\to\infty} \frac{1}{\mu(2Q_k)} \int_{Q_k} |f| \,d\mu >
\frac{\lambda}{\beta_d}.$$
\erem


\subsection{The definition of $\rbmo(\mu)$}

Given a cube $Q\subset\R^d$, let $N$ be the
smallest integer $\geq0$ such that $2^N Q$ is doubling.
We denote this cube by $\wt{Q}$ (recall that this cube $\wt{Q}$ exists
because otherwise the growth condition \rf{creix} on $\mu$ would fail).

Let $\rho>1$ be some fixed constant. We say that $f\in L^1_{loc}(\mu)$ is in
$\rbmo(\mu)$ if there exists some constant
$C_3$ such that for any cube $Q$ (centered at some point of $\supp(\mu)$),
\begin{equation} \label{def1}
\frac{1}{\mu(\rho Q)}\int_Q |f-m_{\wt{Q}}f|\,d\mu \leq  C_3
\end{equation}
and
\begin{equation} \label{def2}
|m_Q f-m_R f|\leq C_3 \,K_{Q,R} \quad \mbox{for any two
doubling cubes $Q\subset R$}.
\end{equation}
The minimal constant $C_3$ is the $\rbmo(\mu)$ norm of $f$ (in fact, it is a
norm in the space of functions modulo additive constants), and it will be
denoted by $\|\cdot\|_*$.

Let us remark that the space $\rbmo(\mu)$ depends on the integer $n$ because
of the definition of the coefficients $K_{Q,R}$.

Notice that if \rf{def1} is satisfied, then \rf{nbmo} also holds. Indeed, for
any cube $Q$ and any $a\in \R$ one has
$$\int_Q |f-m_Q f|\,d\mu \leq 2\int_Q |f-a|\,d\mu.$$
In particular this holds for $a=m_{\wt{Q}} f$. So, the condition \rf{def1}
is stronger than \rf{nbmo}. Moreover, in the definition of $\rbmo(\mu)$
we ask also the regularity condition \rf{def2}.

Observe also that, as a consequence of \rf{def2}, if $Q\subset R$ are doubling
cubes with comparable size, then
\begin{equation} \label{prop1}
|m_Q f-m_R f|\leq C\,\|f\|_*,
\end{equation}
taking into account the property 2 of Lemma \ref{kq}.

\brem \label{rem1.5}
In fact, \rf{prop1}
also holds for any two doubling cubes with comparable sizes such
that $\dist(Q,R)\lessapprox l(Q)$. To see this, let $Q_0$ be a cube concentric
with $Q$, containing $Q$ and $R$, and such that $l(Q_0)\approx l(Q)$. Then
$K(Q_0,\wt{Q_0})\leq C$, and thus we have $K(Q,\wt{Q_0}) \leq C$ and
$K(R,\wt{Q_0})\leq C$ (we have used the properties 1, 2 and 3 of Lemma
\ref{kq}). Then $|m_Q f-m_{\wt{Q_0}} f|\leq C\,\|f\|_*$ and
$|m_R f-m_{\wt{Q_0}} f|\leq C\,\|f\|_*$. So \rf{prop1} holds.
\erem

Let us remark that the condition \rf{prop1} is not satisfied, in general,
by functions of the bigger space $\bmo_\rho(\mu)$ and cubes $Q,\,R$ as above.

We have the following properties:

\begin{propo}  \label{prope}

\begin{enumerate}
\item[1.]$\rbmo(\mu)$ is a Banach space of functions (modulo additive
constants).
\item[2.] $L^\infty(\mu)\subset\rbmo(\mu)$, with $\|f\|_*
\leq 2\|f\|_{L^\infty(\mu)}$.
\item[3.] If $f\in \rbmo(\mu)$, then $|f|\in \rbmo(\mu)$ and $\|\, |f|
\,\|_* \leq C\,\|f\|_*$.
\item[4.] If $f,g\in \rbmo(\mu)$, then $\min(f,g),\, \max(f,g) \in\rbmo(\mu)$
and
$$\|\min(f,g)\|_*,\, \|\max(f,g)\|_* \leq C\,(\|f\|_* + \|g\|_*).$$
\end{enumerate}
\end{propo}

\begin{proof} The properties 1 and 2 are easy to check. The third property
is also easy to prove with the aid of Lemma \ref{normeq2} below. The fourth
property follows from the third.
\end{proof}

Before showing that CZO's which are bounded on $L^2(\mu)$ are also bounded
from
$L^\infty(\mu)$ into $\rbmo(\mu)$, we will see other equivalent norms for
$\rbmo(\mu)$.
Suppose that for a given a function $f\in L^1_{loc}(\mu)$ there exist
some constant $C_4$ and
a collection of numbers $\{f_Q\}_Q$ (i.e. for each cube $Q$, there exists
$f_Q\in\R$) such that
\begin{equation} \label{def1'}
\sup_Q \frac{1}{\mu(\rho Q)} \int_Q |f(x)-f_Q|\,d\mu(x) \leq C_4,
\end{equation}
and,
\begin{equation} \label{def2'}
|f_Q-f_R|\leq C_4\,K_{Q,R} \quad \mbox{for any two cubes $Q\subset R$}.
\end{equation}
Then, we write $\|f\|_{**}=\inf C_4$, where the infimum is taken over all
the constants $C_4$ and all the numbers $\{f_Q\}$ safisfying \rf{def1'} and
\rf{def2'}. It is esily checked that $\|\cdot\|_{**}$ is a norm in the
space of functions modulo constants.

The definition of the norm $\|\cdot\|_{**}$ depends on the constant $\rho$
chosen in \rf{def1'} (the same occurs for $\|\cdot\|_*$). However,
if we write $\|\cdot\|_{**,(\rho)}$ instead of $\|\cdot\|_{**}$, we have

\begin{lemma} \label{normeq1}
The norms $\|\cdot\|_{**,(\rho)}$, $\rho>1$, are equivalent.
\end{lemma}

\begin{proof} Let $\rho>\eta>1$ be some fixed constants.
Obviously, $\|f\|_{**,(\rho)}\leq \|f\|_{**,(\eta)}$. So we only have to show
$\|f\|_{**,(\eta)}\leq C\,\|f\|_{**,(\rho)}$. It is enough to prove that
for a fixed collection of numbers $\{f_Q\}_Q$ satisfying
$$\sup_Q \frac{1}{\mu(\rho Q)} \int_Q |f(x)-f_Q|\,d\mu(x) \leq
2\,\|f\|_{**,(\rho)}$$
and
$$|f_Q-f_R|\leq 2\,K_{Q,R}\,\|f\|_{**,(\rho)}\quad \mbox{for any two cubes
$Q\subset R$,}$$
we have
\begin{equation} \label{nn10}
\frac{1}{\mu(\eta Q_0)}\int_{Q_0} |f-f_{Q_0}|\,d\mu \leq
C\,\|f\|_{**,(\rho)} \quad \mbox{for any fixed cube $Q_0$}.
\end{equation}

For any $x\in Q_0\cap \supp(\mu)$, let $Q_x$ be a cube centered at $x$ with
side length $\frac{\eta-1}{10\rho}\,l(Q_0)$. Then
$l(\rho Q_x) = \frac{\eta-1}{10}\,l(Q_0),$
and so $\rho Q_x\subset \eta Q_0$.
By Besicovich's covering theorem, there exists a family of
points $\{x_i\}_i\subset Q_0\cap \supp(\mu)$ such that the cubes
$\{Q_{x_i}\}_i$ form an almost disjoint covering of $Q_0\cap \supp(\mu)$.
Since $Q_{x_i}$ and $Q_0$ have comparable
sizes,
$$|f_{Q_{x_i}} - f_{Q_0}| \leq C\,\|f\|_{**,(\rho)},$$
with $C$ depending on $\eta$ and $\rho$. Therefore,
\begin{eqnarray*}
\int_{Q_{x_i}} |f-f_{Q_0}|\,d\mu & \leq &
\int_{Q_{x_i}} |f-f_{Q_{x_i}}|\,d\mu +  |f_{Q_0}-f_{Q_{x_i}}|\,\mu(Q_{x_i})\\
& \leq & C\,\|f\|_{**,(\rho)}\,\mu(\rho Q_{x_i}).
\end{eqnarray*}
Then we get
$$\int_{Q_0} |f-f_{Q_0}|\,d\mu \leq
\sum_i \int_{Q_{x_i}} |f-f_{Q_0}|\,d\mu \leq
C\,\|f\|_{**,(\rho)} \sum_i \mu(\rho Q_{x_i}).$$
Since $\rho Q_{x_i}\subset \eta Q_0$ for all $i$, we obtain
$$\int_{Q_0} |f-f_{Q_0}|\,d\mu \leq
C\,\|f\|_{**,(\rho)}\, \mu(\eta Q_0)\,N,$$
where $N$ is the number of cubes of the Besicovich covering. Now it is easy
to check that $N$ is bounded some constant depending only on $\eta$,
$\rho$ and $d$: If $\LL^d$ is the Lebesgue measure on $\R^d$ and
$B_d$ is the Besicovich constant in $\R^d$, we have
$$N\, \LL^d(Q_{x_i}) = \sum_i \LL^d(Q_{x_i})
\leq B_d\,\LL^d(\eta Q_0).$$
Thus
$$N\leq \frac{B_d\,\LL^d(\eta Q_0)}{\LL^d(Q_{x_i})} =
B_d\,\left(\frac{10\eta\rho}{\eta-1}\right)^d,$$
and \rf{nn10} holds. \end{proof}

\brem \label{rem22}
In fact, in the preceeding lemma we have seen that
if $C_f$ is some constant and $\{f_Q\}_Q$ is some fixed collection of numbers
satisfying
$$\sup_Q \frac{1}{\mu(\rho Q)} \int_Q |f(x)-f_Q|\,d\mu(x) \leq C_f$$
and
$$|f_Q-f_R|\leq K_{Q,R}\,C_f\quad \mbox{for any two cubes
$Q\subset R$,}$$
then for the same numbers $\{f_Q\}_Q$, for any $\eta>1$ we have
$$\sup_Q \frac{1}{\mu(\eta Q)} \int_Q |f(x)-f_Q|\,d\mu(x) \leq
C\,C_f,$$
with $C$ depending on $\eta$.
\erem

We also have:

\begin{lemma}  \label{normeq2}
For a fixed $\rho>1$, the norms $\|\cdot\|_{*}$  and $\|\cdot\|_{**}$ are
equivalent.
\end{lemma}

\begin{proof} Let $f\in L^1_{loc}(\mu)$.
To see that $\|f\|_{**}\leq C\,\|f\|_{*}$
we set $f_Q=m_\wt{Q} f$ for all cubes $Q$. Then \rf{def1'}
holds with $C_4=\|f\|_*$. Let us check that the second condition \rf{def2'} is
also satisfied. We have to prove that for
any two cubes $Q\subset R$,
\begin{equation} \label{nn1}
|m_{\wt{Q}}f - m_{\wt{R}}f| \leq C\,K_{Q,R}\,\|f\|_*.
\end{equation}
Notice that if $\wt{Q}\subset\wt{R}$, then
$$|m_{\wt{Q}}f - m_{\wt{R}}f| \leq K_{\wt{Q},\wt{R}}\,\|f\|_*,$$
because $\wt{Q},\wt{R}$ are doubling. So \rf{nn1} follows if
$K_{\wt{Q},\wt{R}}\leq C\,K_{Q,R}$.
However, in general, $Q\subset R$ does not imply $\wt{Q}\subset\wt{R}$, and
so we have to modify the argument.

Suppose first that $l(\wt{R}) \geq l(\wt{Q})$. Then $\wt{Q} \subset 4 \wt{R}$.
We denote ${R_0} = \wt{4\wt{R}}$. Then we have
\begin{equation}  \label{nn2}
|m_\wt{Q} f - m_\wt{R} f| \leq |m_\wt{Q} f - m_{R_0} f| +
|m_{R_0} f - m_\wt{R} f|.
\end{equation}
Using the properties of Lemma \ref{kq} repeatedly, we get
\begin{eqnarray*}
K_{\wt{Q},{R_0}} & \leq & C\,K_{Q,{R_0}} \leq C\,(K_{Q,R}+K_{R,{R_0}})\\
& \leq& C\,(K_{Q,R}+K_{R,\wt{R}}+ K_{\wt{R},4\wt{R}} + K_{4\wt{R},{R_0}})
\leq C\,K_{Q,R}.
\end{eqnarray*}
Since $\wt{Q}\subset{R_0}$ and they are doubling cubes, we have
$$|m_\wt{Q} f - m_{R_0} f| \leq K_{\wt{Q},{R_0}}\,\|f\|_* \leq
C\,K_{Q,R}\,\|f\|_*.$$
Now we are left with the second term on the right hand side of
\rf{nn2}. We have
$$K_{\wt{R},{R_0}} \leq C(K_{\wt{R},4\wt{R}}+ K_{4\wt{R},{R_0}}) \leq C
\leq C\,K_{Q,R}.$$
Due to the fact that $\wt{R}\subset{R_0}$ are doubling cubes,
$$|m_{R_0} f - m_\wt{R} f| \leq K_{\wt{R},{R_0}}\,\|f\|_* \leq
C\,K_{Q,R} \,\|f\|_*,$$
and, by \rf{nn2}, we get that \rf{nn1} holds in this case.

Assume now $l(\wt{R}) < l(\wt{Q})$. Then $\wt{R}\subset4\wt{Q}$.
There exists some $m\geq1$ such that
$l(\wt{R})\geq l(2^mQ)/10$ and $\wt{R}\subset 2^mQ \subset 4\wt{Q}$.
Since $\wt{R}$ and $2^mQ$ have comparable sizes, we have $K_{\wt{R},2^mQ}\leq
C$. Then, if we denote ${Q_0} = \wt{4\wt{Q}}$, we get
$$K_{\wt{R},{Q_0}} \leq C\,(K_{\wt{R},2^mQ} + K_{2^mQ,4\wt{Q}} +
K_{4\wt{Q},{Q_0}}) \leq C.$$
Also, $$K_{\wt{Q},{Q_0}} \leq C\,(K_{\wt{Q},4\wt{Q}} +
K_{4\wt{Q},{Q_0}})\leq C.$$
Therefore,
\begin{eqnarray*}
|m_\wt{Q} f - m_\wt{R} f| & \leq &
|m_\wt{Q} f - m_{Q_0} f| +
|m_{Q_0} f - m_\wt{R} f| \\
& \leq & K_{\wt{Q},{Q_0}}\,\|f\|_*
+ K_{\wt{R},{Q_0}}\,\|f\|_* \leq C\,\|f\|_*
\leq C\,K_{Q,R}\,\|f\|_*.
\end{eqnarray*}

Now we have to check that $\|f\|_*\leq C\,\|f\|_{**}$. If $Q$ is a doubling
cube, since \rf{def1'} holds with $\rho=2$ (by Lemma \ref{normeq1}), we have
$$|f_Q-m_Qf| = \left|\frac{1}{\mu(Q)} \int_Q (f-f_Q)\,d\mu\right| \leq
\|f\|_{**}\,\frac{\mu(2Q)}{\mu(Q)} \leq C\,\|f\|_{**}.$$
Therefore, for any cube $Q$ (non doubling, in general), using
$K_{Q,\wt{Q}}\leq C$ we get
$$|f_Q-m_{\wt{Q}}f|\leq |f_Q-f_{\wt{Q}}| + |f_{\wt{Q}}-m_{\wt{Q}}f|
\leq C \,\|f\|_{**}.$$
Thus
\begin{eqnarray*}
\frac{1}{\mu(\rho Q)} \int_Q |f(x)-m_{\wt{Q}}f|\,d\mu(x) & \leq &
\frac{1}{\mu(\rho Q)} \int_Q |f(x)-f_Q|\,d\mu(x) \\ &&\mbox{} +
\frac{1}{\mu(\rho Q)} \int_Q |f_Q-m_{\wt{Q}}f|\,d\mu(x) \\
& \leq & C\,\|f\|_{**}.
\end{eqnarray*}
It only remains to show that \rf{def2} also holds with $C\,\|f\|_{**}$ instead
of $C_3$. This follows easily. Indeed, if $Q\subset R$ are
doubling cubes, we have
\begin{eqnarray*}
|m_Q f-m_R f| & \leq & |m_Q f-f_Q| + |f_Q - f_R| + |f_R-m_R f| \\
& \leq &
C\, \|f\|_{**} + K_{Q,R}\,\|f\|_{**} \leq C\,K_{Q,R}\,\|f\|_{**}.
\end{eqnarray*}
\end{proof}

\brem  \label{remnou}
By the preceeding lemma, it is easily seen
that we obtain equivalent norms and the same space $\rbmo(\mu)$ if we replace
$(2,\beta_d)$-doubling cubes in the definition
of the space $\rbmo(\mu)$ by $(\alpha,\beta)$-doubling cubes, for any choice
of $\alpha>1$ and $\beta>\alpha^n$. We have taken $(2,\beta_d)$-doubling
cubes in the definition of $\|\cdot\|_*$ (and not $(2,2^{n+1})$-doubling, say)
because to prove some of the results below it will be necessary to work
with doubling cubes having the properties explained in Remark \ref{rem9}.

On the other hand, by Lemmas \ref{normeq1} and \ref{normeq2},
the definition of $\rbmo(\mu)$ does not depend on the
number $\rho>1$ chosen in \rf{def1}. So, throughout the rest of the paper
we will assume that the constant $\rho$ in the definition of $\rbmo(\mu)$ is
$2$.

Also, it can be seen that we also obtain equivalent definitions 
for the space $\rbmo(\mu)$ if instead of cubes centered at points in
$\supp(\mu)$, we consider all the cubes in $\R^d$ (with sides parallel to the
axes). Furthermore, it does not matter if we take balls instead of cubes. 

Notice that in Lemma \ref{normeq2} we have shown that if we choose
$f_Q=m_{\wt{Q}}f$ for all cubes $Q$, then \rf{def1'} and \rf{def2'} are
satisfied with
$C_4 = C\,\|f\|_*$.
\erem

Other possible ways of defining $\rbmo(\mu)$ are shown in the following lemma.

\begin{lemma} \label{fifi}
Let $\rho>1$ be fixed. For a function $f\in L^1_{loc}(\mu)$, the following are
equivalent:

\begin{itemize}
\item[a)] $f\in\rbmo(\mu)$.

\item[b)] There exists some constant $C_b$ such that for any cube $Q$
\begin{equation} \label{def1''}
\int_Q |f-m_Q f|\,d\mu \leq  C_b\, \mu(\rho Q)
\end{equation}
and
\begin{equation} \label{def2''}
|m_Q f-m_R f|\leq C_b \,K_{Q,R} \left(\frac{\mu(\rho Q)}{\mu(Q)} +
\frac{\mu(\rho R)}{\mu(R)}\right) 
\quad \mbox{for any two cubes $Q\subset R$}.
\end{equation}

\item[c)] There exists some constant $C_c$ such that for any {\em doubling} cube $Q$
\begin{equation} \label{def1'''}
\int_Q |f-m_Q f|\,d\mu \leq  C_c\, \mu(Q)
\end{equation}
and
\begin{equation} \label{def2'''}
|m_Q f-m_R f|\leq C_c \,K_{Q,R}  
\quad \mbox{for any two {\em doubling} cubes $Q\subset R$}.
\end{equation}
\end{itemize}

Moreover, the best constants $C_b$ and $C_c$ are comparable to the
$\rbmo(\mu)$ norm of $f$. 
\end{lemma}

\begin{proof}
Assume $\rho=2$ for simplicity. 
First we show a) $\Rightarrow$ b).
If $f\in\rbmo(\mu)$, then \rf{def1''} holds with
$C_b=2\|f\|_*$. Moreover, for any cube $Q$ we have
\begin{equation}  \label{fifi1}
|m_Q f - m_{\wt{Q}}f|\leq m_Q (|f - m_{\wt{Q}}f|) \leq \|f\|_*
\frac{\mu(2Q)}{\mu(Q)}.
\end{equation}
Therefore,
$$|m_Q f - m_Rf|\leq |m_Q f - m_\wt{Q}f| + |m_\wt{Q} f - m_\wt{R}f|+
|m_R f - m_\wt{R}f|.$$
The second term on the right hand side is estimated as \rf{nn1} in the
preceeding lemma. For the first and third terms on the right, we apply
\rf{fifi1}. So we get,
\begin{eqnarray*}
|m_Q f - m_Rf| & \leq & \left(C\,K_{Q,R} + \frac{\mu(2Q)}{\mu(Q)} +
\frac{\mu(2R)}{\mu(R)} \right)\, \|f\|_*\\
& \leq  &
C\,K_{Q,R} \left( \frac{\mu(2Q)}{\mu(Q)} +
\frac{\mu(2R)}{\mu(R)} \right)\, \|f\|_*.
\end{eqnarray*}
Thus $f$ satisfies \rf{def2''} too.

The implication b) $\Rightarrow$ c) is easier: One only has to consider
doubling cubes in b).

Let us see now c)$\Rightarrow$ a).
Let $Q$ be some cube, non doubling in general. We only have to show that
\rf{def1} holds. We know that for $\mu$-almost all $x\in Q$ there exists
some {\em doubling} cube centered at $x$ with sidelength $2^{-k}\,l(Q)$,
for some $k\geq1$. We denote by $Q_x$ the biggest cube satisfying these
properties. Observe that $K_{Q_{x},\wt{Q}}\leq C$, and then
\begin{equation}  \label{tyu}
|m_{Q_x} f - m_{\wt{Q}} f|\leq C\,C_c.
\end{equation}

By Besicovich's covering theorem, there are points $x_i\in Q$ such that
$\mu$-almost all $Q$ is covered by a family of cubes 
$\{Q_{x_i}\}_i$ with bounded overlap. By \rf{tyu}, using that 
$Q_{x_i}\subset 2Q$, we get
\begin{eqnarray*}
\int_Q|f- m_\wt{Q} f|\,d\mu &\leq & \sum_i
\int_{Q_{x_i}} |f- m_\wt{Q} f|\,d\mu \\
&\leq & \sum_i \int_{Q_{x_i}} |f- m_{Q_{x_i}} f|\,d\mu +  
\sum_i |m_\wt{Q}f- m_{Q_{x_i}} f|\, \mu(Q_{x_i}) \\
& \leq & C\, C_c\, \mu(2Q).
\end{eqnarray*}
\end{proof}


\subsection{Boundedness of CZO's from $L^\infty(\mu)$ into $\rbmo(\mu)$}

Now we are going to see that if a CZO is bounded on $L^2(\mu)$, then
it is bounded from $L^\infty(\mu)$ into $\rbmo(\mu)$. In fact, we will
replace the assumption of $L^2(\mu)$ boundedness by another weaker assumption.

\begin{theorem} \label{aco1}
If for any cube $Q$ and any function $a$ supported on $Q$
\begin{equation}  \label{hip1}
\int_Q |T_\ve a|\,d\mu \leq C\,\|a\|_{L^\infty}\,\mu(\rho Q)
\end{equation}
uniformly on $\ve>0$, then $T$ is bounded from $L^\infty(\mu)$ into
$\rbmo(\mu)$.
\end{theorem}

Let us remark that when we say that $T$ is bounded from $L^\infty(\mu)$ into
$\rbmo(\mu)$, we mean that the operators $T_\ve$, $\ve>0$, are uniformly
bounded from $L^\infty(\mu)$ into $\rbmo(\mu)$.

\vspace{3mm}
\begin{proof}
First we will show that if $f\in L^\infty(\mu)\cap L^{p_0}(\mu)$ for
some $p_0\in [1,\infty)$, then
\begin{equation} \label{mnm}
\|T_\ve f\|_{\rbmo(\mu)} \leq C\|f\|_{L^\infty(\mu)}.
\end{equation}
We will use the characterization of $\rbmo(\mu)$ given by \rf{def1''} and \rf{def2''}
in Lemma \ref{fifi}.

The condition \rf{def1''} follows by standard methods. The same proof that
shows that $T_\ve f\in \bmo(\mu)$ when $\mu$ is a doubling measure works here.
We omit the details.

Let us see how \rf{def2''} follows. For simplicity, we assume $\rho=2$. We have 
to show that if $Q\subset R$, then
$$|m_Q(T_\ve f) - m_R(T_\ve f)|\leq C\,K_{Q,R}\,
\left(\frac{\mu(2Q)}{\mu(Q)} +
\frac{\mu(2R)}{\mu(R)}\right)\,\|f\|_{L^\infty(\mu)}.$$
Recall that $N_{Q,R}$ is the first integer $k$ such that $2^{k}Q\supset
R$. We denote $Q_R = 2^{N_{Q,R}+1}Q$. Then, for $x\in Q$ and $y\in R$, we
set
\begin{eqnarray*}  
T_\ve f (x) - T_\ve f(y) & = & T_\ve f \,\chi_{2Q}(x) +\sum_{k=1}^{N_{Q,R}}
T_\ve f \,\chi_{2^{k+1}Q \setminus 2^{k}Q}(x) + T_\ve f
\,\chi_{\R^d\setminus Q_R}(x) \\
&& \mbox{} - \left(
T_\ve f \,\chi_{Q_R}(y) + T_\ve f
\,\chi_{\R^d\setminus Q_R}(y) \right).
\end{eqnarray*}
Since
$$|T_\ve f \,\chi_{\R^d\setminus Q_R}(x) -
T_\ve f \,\chi_{\R^d\setminus Q_R}(y)| \leq C\,\|f\|_{L^\infty(\mu)},$$
we get
\begin{eqnarray} \label{eq4}
|T_\ve f (x) - T_\ve f(y)| & \leq & |T_\ve f \,\chi_{2Q}(x)| +
C\,\sum_{k=1}^{N_{Q,R}} \frac{\mu(2^{k+1}Q)}{
l(2^{k+1}Q)^n}\,\|f\|_{L^\infty(\mu)} \nonumber \\
&& \mbox{} + |T_\ve f \,\chi_{Q_R}(y)|
+ C\,\|f\|_{L^\infty(\mu)}.
\end{eqnarray}
Now we take the mean over $Q$ for $x$, and over $R$ for $y$.
Using the $L^2(\mu)$ boundedness of $T$, we obtain
\begin{eqnarray*}
m_Q(|T_\ve f \,\chi_{2Q}|) & \leq & \left(\frac{1}{\mu(Q)} \int_Q 
|T_\ve f \,\chi_{2Q}|^2\,\, d\mu\right)^{1/2} \\
& \leq & C\,\left(\frac{\mu(2Q)}{\mu(Q)}\right)^{1/2} \|f\|_{L^\infty(\mu)} \\
& \leq & C\,\frac{\mu(2Q)}{\mu(Q)} \|f\|_{L^\infty(\mu)}.
\end{eqnarray*}
For $R$ we write
$$m_R(|T_\ve f \,\chi_{Q_R}|) \leq
m_R(|T_\ve f \,\chi_{Q_R\cap 2R}|)  +  m_R(|T_\ve f \,\chi_{Q_R\setminus 2R}|).$$
The estimate for the first term on the right hand side is similar to the
previous estimate for $Q$:
\begin{eqnarray*}
m_R(|T_\ve f \,\chi_{Q_R\cap 2R}|) & \leq &
\left(\frac{1}{\mu(R)} \int_R
|T_\ve f \,\chi_{Q_R\cap 2R}|^2\, d\mu\right)^{1/2} \\
& \leq & C\,\left(\frac{\mu(Q_R\cap 2R)}{\mu(R)}\right)^{1/2}
\|f\|_{L^\infty(\mu)}\\
& \leq & C\,\frac{\mu(2R)}{\mu(R)} \|f\|_{L^\infty(\mu)}.
\end{eqnarray*}
On the other hand, since $l(Q_R)\approx l(R)$, we have
$m_R(|T_\ve f \,\chi_{Q_R\setminus 2R}|)\leq C\, \|f\|_{L^\infty(\mu)}$.
Therefore,
\begin{eqnarray*}
|m_Q(T_\ve f) - m_R(T_\ve f)| & \leq & C\,
\sum_{k=1}^{N_{Q,R}} \frac{\mu(2^{k+1}Q)}{
l(2^{k+1}Q)^n}\,\|f\|_{L^\infty(\mu)} \\
&& \mbox{} +
C\,\left(\frac{\mu(2Q)}{\mu(Q)} + \frac{\mu(2R)}{\mu(R)}\right)
\,\|f\|_{L^\infty(\mu)} \\
& \leq & C\,K_{Q,R}\,\left(\frac{\mu(2Q)}{\mu(Q)} +
\frac{\mu(2R)}{\mu(R)}\right) \,\|f\|_{L^\infty(\mu)}.
\end{eqnarray*}
So we have proved that \rf{mnm} holds for $f\in L^\infty(\mu)\cap
L^{p_0}(\mu)$.

\vv
If $f\not\in L^p(\mu)$ for all $p\in[1,\infty)$, then the integral
$\int_{|x-y|>\ve} k(x,y)\,f(y)\,d\mu(y)$ may be not convergent.
The operator $T_\ve$ can be extended to the whole space $L^\infty(\mu)$
following the usual arguments:
Given a cube $Q_0$ centered at the origin with side length $>3\ve$, we write
$f=f_1+f_2$, with $f_1=f\,\chi_{2Q_0}$. For $x\in Q_0$, we define
$$T_\ve f(x) = T_\ve f_1(x) + \int (k(x,y)-k(0,y))\,f_2(y)\,d\mu(y).$$
Now both integrals in this equation are convergent and
with this definition one can check that \rf{mnm} holds too, with
arguments similar to the case
$f\in L^\infty(\mu)\cap L^{p_0}(\mu)$.
\end{proof}

Let us remark that in Theorem \ref{equi}
we will see that the condition \rf{hip1} holds
if and only if $T$ is bounded from $L^\infty(\mu)$ into $\rbmo(\mu)$.


\subsection{Examples}

\bexam \label{ex1} Assume $d=2$ and $n=1$. So we can think that we are in
the complex plane and $T$ is the Cauchy transform. Let $E\subset \C$ be a
$1$-dimensional Ahlfors-David (AD) regular set. That is,
$$C^{-1}\,r\leq \HH^1(E\cap B(x,r)) \leq C\,r\quad \mbox{for all $x\in E$,
$0<r\leq \diam(E)$.}$$
(Here $\HH^1$ stands for the $1$-dimensional Hausdorff measure.)
We set $\mu=\HH^1_{\mid E}$. Notice that $\mu$ is a doubling measure.

For any $Q$ centered at
some point of $\supp(\mu)$, one has
$$\mu(2^k Q)\approx l(2^k Q)$$
if $l(2^k Q)\leq \diam(E)$. Then, given $Q\subset R$,
it is easy to check that if $l(R)\leq \diam(E)$,
\begin{equation}  \label{eqk1}
K_{Q,R} \approx 1 + \log\frac{l(R)}{l(Q)},
\end{equation}
and if $l(R)> \diam(E)$,
\begin{equation}  \label{eqk2}
K_{Q,R} \approx 1 + \log\frac{\diam(E)}{l(Q)}.
\end{equation}
So, in this case, we have $\rbmo(\mu) = \bmo(\mu)$, since any function
$f\in \bmo(\mu)$ satisfies \rf{def1'} and \rf{def2'}, with $f_Q=m_Q f$ for all
cubes $Q$. Notice that \rf{def2'} holds
because of \rf{eqk1} and \rf{eqk2}.
\eexam

\bexam \label{ex2} We assume again $d=2$ and $n=1$. Let $\mu$ be the planar
Lebesgue measure restricted to the unit square $[0,1]\times [0,1]$. This
measure is doubling, but not AD-regular (for $n=1$). Now one can check
that the coefficients $K_{Q,R}$ are uniformly bounded. That is, for any two
squares $Q\subset R$,
$$K_{Q,R} \approx 1.$$
Let us take $R_0=[0,1]^2$ and $Q\subset R_0$. Then, if $f\in \rbmo(\mu)$,
$$|m_Q (f - m_{R_0} f)| = |m_Q f - m_{R_0} f|\leq K_{Q,R_0}\,\|f\|_*
\leq C\,\|f\|_*.$$
Since this holds for any square $Q\subset R_0$, by the Lebesgue
differentiation theorem $f - m_{R_0} f$ is a bounded function, with
$$\|f - m_{R_0} f\|_{L^\infty(\mu)} \leq \|f\|_*.$$
Therefore, now $\rbmo(\mu)$ coincides with $L^\infty(\mu)$ modulo constants
functions, which is strictly smaller than $\bmo(\mu)$.
\eexam

\bexam \label{ex3} This example is borrowed in part from
\cite{NTV*}. Suppose $d=2$ (i.e. we are in the complex plane) and $n=1$.
Let $\mu$ be a measure on the real axis such that in the intervals $[-2,-1]$
and $[1,2]$ is the linear Lebesgue measure, on the interval $[-1/2,1/2]$ is
the linear Lebesgue measure times $\ve$, with $\ve>0$ very small, and $\mu=0$
elsewhere.
We consider the function $f =
\ve^{-1}\,(\chi_{[1/4,1/2]}-\chi_{[-1/2,-1/4]})$.
It is easily checked that for $\rho\leq 2$,
$$\|f\|_{\bmo_\rho(\mu)} \approx \ve^{-1},$$
while for $\rho=5$,
$$\|f\|_{\bmo_5(\mu)} \approx 1.$$
On the other hand, the $\rbmo(\mu)$ norm of $f$ is
$$\|f\|_* \approx \ve^{-1},$$
since
$$C\,\|f\|_* \geq |m_{[-2,2]}f - m_{[1/4,1/2]} f| = \ve^{-1}$$
and $\|f\|_*\leq C\,\|f\|_{L^\infty(\mu)} \leq C\,\ve^{-1}$.
\eexam


\section{The inequality of John-Nirenberg}

The following result is a version of John-Nirenberg's inequality suitable
for the space $\rbmo(\mu)$. To prove it we will adapt the arguments in
\cite[p.31-32]{Journe} to the present situation.

\begin{theorem} \label{JN}
Let $f\in \rbmo(\mu)$ and let $\{f_Q\}_Q$ be a collection of numbers
satisfying
\begin{equation}  \label{mn3}
\sup_Q \frac{1}{\mu(2 Q)} \int_Q |f(x)-f_Q|\,d\mu(x) \leq C\,\|f\|_*
\end{equation}
and
\begin{equation}  \label{mn4}
|f_Q-f_R|\leq C\,K_{Q,R}\,\|f\|_* \quad \mbox{for any two cubes $Q\subset
R$.}
\end{equation}
Then, for any cube $Q$ and any $\lambda >0$ we have
\begin{equation}  \label{mn88}
\mu\{x\in Q:|f(x)-f_Q|>\lambda\} \leq C_5\,\mu(\rho
Q)\,\exp\left({\frac{-C_6\,\lambda}{\|f\|_*}}\right),
\end{equation}
with $C_5$ and $C_6$ depending on the constant $\rho>1$ (but not on
$\lambda$).
\end{theorem}

To prove this theorem, we will use the following straightforward result:

\begin{lemma} \label{facil}
Let $f\in \rbmo(\mu)$ and let $\{f_Q\}_Q$ be a collection of numbers
satisfying \rf{mn3} and \rf{mn4}. If $Q$ and $R$ are cubes
such that $l(Q)\approx l(R)$ and $\dist(Q,R) \lessapprox l(Q)$, then
$$|f_Q-f_R|\leq C\,\|f\|_*.$$
\end{lemma}

\begin{proof} Let $R'$ be the smallest cube concentric with $R$ containing $Q$ and
$R$. Since $l(Q)\approx l(R') \approx l(R)$, we have $K_{Q,R'}\leq C$ and
$K_{R,R'}\leq C$. Then,
$$|f_Q-f_R| \leq |f_Q-f_{R'}| +  |f_R-f_{R'}| \leq
C\,(K_{Q,R} + K_{R,R'}) \,\|f\|_* \leq C\,\|f\|_*.$$
\end{proof}

We will use the following lemma too.

\begin{lemma}  \label{lemaq}
Let $f\in\rbmo(\mu)$. Given $q>0$, we set
$$f_q(x) = \left\{ \begin{array}{ll}
f(x) \hspace{3mm} & \mbox{if $|f(x)| \leq q$,} \\
& \\
\ds q\, \frac{f(x)}{|f(x)|} \hspace{3mm} & \mbox{if $|f(x)|>q$.}
\end{array}
\right. $$
Then $f_q \in \rbmo(\mu)$, with $\|f_q\|_* \leq C\,\|f\|_*$.
\end{lemma}

\begin{proof}
For any function $g$, we set $g=g_+-g_-$, with $g_+ = \max(g,0)$ and
$g_-= -\min(g,0)$.

By Proposition \ref{prope}, $\|f_+\|_*,\,\|f_-\|_*\leq C\,\|f\|_*$.
Since $f_{q,+} = \min(f_+,q)$ and $f_{q,-} = \min(f_-,q)$, we have
$\|f_{q,+}\|_*, \|f_{q,-}\|_* \leq C\,\|f\|_*$. Thus $\|f_q\|_* \leq
\|f_{q,+}\|_* + \|f_{q,-}\|_* \leq C\,\|f\|_*$.
\end{proof}

\brem \label{remqqq}
Let $f\in \rbmo(\mu)$ and let $\{f_Q\}_Q$ be a collection of numbers
satisfying \rf{mn3} and \rf{mn4}. We set $f_{Q,+} = \max(f_Q,0)$ and 
$f_{Q,-} = -\min(f_Q,0)$ and we take
$$f_{q,Q} = \min(f_{Q,+},q) - \min(f_{Q,-},q).$$
It is easily seen that
$$\sup_Q \frac{1}{\mu(2 Q)} \int_Q |f_q(x)-f_{q,Q}|\,d\mu(x) \leq C\,\|f\|_*$$
and
$$
|f_{q,Q}-f_{q,R}|\leq C\,K_{Q,R}\,\|f\|_* \quad \mbox{for any two cubes $Q\subset
R$.}$$
\erem

\vspace{5mm}
\begin{proof}[Proof of Theorem \ref{JN}.]
We will prove \rf{mn88} for $\rho=2$. The proof for other values of $\rho$ is
similar. Recall that if \rf{mn3} and \rf{mn4} are satisfied, then \rf{mn3}
is also satisfied subsituting ``$\mu(2Q)$'' by ``$\mu(\rho Q)$'', for any
$\rho>1$ (see Remark \ref{rem22}).

Let $f\in\rbmo(\mu)$. Assume first that $f$ is bounded.
Let $Q_0$ be some fixed cube in $\R^d$. We write $Q_0'=\frac{3}{2}Q_0$.

Let $B$ be some positive constant which will be fixed later.
By Remark \ref{rem9}, for $\mu$-almost any $x\in Q_0$ such that
$|f(x)-f_{Q_0}|>B\,\|f\|_*$,
there exists some doubling
cube $Q_x$ centered at $x$ satisfying
\begin{equation}  \label{mn1}
m_{Q_x}(|f-f_{Q_0}|)>B\,\|f\|_*.
\end{equation}
Moreover, we may assume that $Q_x$ is the biggest
doubling cube satisfying \rf{mn1} with
side length $2^{-k}\,l(Q_0)$ for some integer $k\geq0$,
with
$$l(Q_x) \leq \frac{1}{10}\,l(Q_0).$$
By Besicovich's covering theorem, there exists an almost 
disjoint subfamily
$\{Q_i\}_i$ of the cubes $\{Q_x\}_x$ such that
\begin{equation}  \label{mn77}
\{x:\, |f(x)-f_{Q_0}|>B\,\|f\|_*\}\subset \bigcup_i Q_i.
\end{equation}
Then, since $Q_i\subset Q_0'$ and
$|f_{Q_0}-f_{Q_0'}|\leq C\,\|f\|_*$, we have
\begin{eqnarray} \label{mn5}
\sum_i \mu(Q_i) &\leq & \sum_i \frac{1}{B\,\|f\|_*}\,
\int_{Q_i} |f-f_{Q_0}|\,d\mu \nonumber\\
& \leq & \frac{C}{B\,\|f\|_*}\,
\int_{Q_0'} |f-f_{Q_0}|\,d\mu \nonumber \\
& \leq & \frac{C}{B\,\|f\|_*}\,|f_{Q_0}-f_{Q_0'}|\,\mu(Q_0') +
\frac{C}{B\,\|f\|_*}\, \int_{Q_0'} |f-f_{Q_0'}|\,d\mu.
\end{eqnarray}
Since  \rf{mn3} is satisfied if we change
``$\mu(2Q)$'' by ``$\mu(\frac{4}{3} Q)$'', we have
$$\int_{Q_0'} |f-f_{Q_0'}|\,d\mu \leq C\,\mu(2Q_0)\,\|f\|_*,$$
and, by \rf{mn5},
$$\sum_i \mu(Q_i) \leq \frac{C\,\mu(2 Q_0)}{B}.$$
So if we choose $B$ big enough,
\begin{equation} \label{mn6}
\sum_i \mu(Q_i) \leq \frac{\mu(2Q_0)}{2\beta_d}.
\end{equation}

Now we want to see that for each $i$ we have
\begin{equation} \label{mn7}
|f_{Q_i}-f_{Q_0}|\leq C_7\,\|f\|_*.
\end{equation}
We consider the cube $\wt{2Q_i}$. If $l(\wt{2Q_i})>10l(Q_0)$, then there
exists
some cube $2^mQ_i$, $m\geq1$, containing $Q_0$ and such that
$l(Q_0) \approx l(2^mQ_i) \leq l(\wt{2Q_i})$. Thus
$$|f_{Q_i}-f_{Q_0}| \leq |f_{Q_i}-f_{2Q_i}| + |f_{2Q_i}-f_{2^mQ_i}|
+ |f_{2^mQ_i}-f_{Q_0}|.$$
The first and third sums on the right hand side are bounded by $C\,\|f\|_*$
because  $Q_i$ and $2Q_i$ on the one hand and  $2^mQ_i$ and $Q_0$ on the
other hand have comparable sizes. The second sum is also bounded by
$C\,\|f\|_*$ due to the fact that there are no
doubling cubes of the form
$2^kQ_i$ between $Q_i$ and $2^mQ_i$, and then $K_{Q_i,2^mQ_i}\leq C$.

Assume now $\frac{1}{10}\,l(Q_0)< l(\wt{2Q_i})\leq10l(Q_0)$. Then
$$|f_{Q_i}-f_{Q_0}| \leq |f_{Q_i}-f_{\wt{2Q_i}}| + |f_\wt{2Q_i}-f_{Q_0}|.$$
Since $\wt{2Q_i}$ and $Q_0$ have comparable sizes, by Lemma \ref{facil}
we have
$|f_\wt{2Q_i}-f_{Q_0}|\leq C\,\|f\|_*$. And since $K_{Q_i,\wt{2Q_i}}
\leq C(K_{Q_i,2Q_i} + K_{2Q_i,\wt{2Q_i}})$, we also have
$|f_{Q_i}-f_{\wt{2Q_i}}|\leq C\,\|f\|_*$. So \rf{mn7} holds in this case too.

If $l(\wt{2Q_i})\leq \frac{1}{10}\,l(Q_0)$, then, by the choice of
$Q_i$, we have
$m_{\wt{2Q_i}}(|f-f_{Q_0}|)\leq B\,\|f\|_*,$
which implies
\begin{equation}  \label{mn23}
|m_{\wt{2Q_i}}(f-f_{Q_0})|\leq B\,\|f\|_*.
\end{equation}
Thus
$$|f_{Q_i} - f_{Q_0}| \leq |f_{Q_i} - f_{\wt{2Q_i}}| +
|f_{\wt{2Q_i}} - m_{\wt{2Q_i}}f| + |m_{\wt{2Q_i}}f - f_{Q_0}|.$$
As above, the term $|f_{Q_i} - f_{\wt{2Q_i}}|$ is bounded by $C\,\|f\|_*$.
The last one equals $|m_{\wt{2Q_i}}(f - f_{Q_0})|$, which is estimated in
\rf{mn23}. For the second one, since $\wt{2Q_i}$ is doubling, we
have
$$|f_{\wt{2Q_i}} - m_{\wt{2Q_i}}f| \leq \frac{1}{\mu(\wt{2Q_i})}\,
\int_\wt{2Q_i} |f-f_\wt{2Q_i}|\,d\mu \leq C\,\frac{\mu(2\cdot\wt{2Q_i})}{
\mu(\wt{2Q_i})}\,\|f\|_*\leq C\,\|f\|_*.$$
So \rf{mn7} holds in any case.

Now we consider the function
$$X(t) = \sup_Q \frac{1}{\mu(2 Q)}\int_Q \exp\left(
|f-f_Q|\,\frac{t}{\|f\|_*} \right)\, d\mu.$$
Since we are assuming that $f$ is bounded, $X(t)<\infty$. By \rf{mn77} and
\rf{mn7} we have
\begin{multline*}
\frac{1}{\mu(2 Q_0)}\int_{Q_0} \exp\left(
|f-f_{Q_0}|\,\frac{t}{\|f\|_*} \right)\, d\mu  \\
\begin{split}
& \leq
\frac{1}{\mu(2 Q_0)}\int_{Q_0\setminus \bigcup_i Q_i} \exp(B\,t) \,d\mu \\
& \mbox{} + \frac{1}{\mu(2 Q_0)} \sum_i
\int_{Q_i} \exp\left(
|f-f_{Q_i}|\,\frac{t}{\|f\|_*} \right)\, d\mu \cdot \exp(C_7\,t) \\
& \leq \exp(B\,t) +  \frac{1}{\mu(2 Q_0)} \sum_i
\mu(2 Q_i)\,X(t)\,\exp(C_7\,t).
\end{split}
\end{multline*}
By \rf{mn6} and taking into account that $\mu(2 Q_i)/\mu(Q_i)\leq \beta_d$,
we get
$$\frac{1}{\mu(2 Q_0)}\int_{Q_0} \exp\left(
|f-f_{Q_0}|\,\frac{t}{\|f\|_*} \right)\, d\mu
\leq \exp(B\,t) + \frac{1}{2}\,X(t)\,\exp(C_7\,t).$$
Thus
$$X(t) \left(1 - \frac{1}{2}\, \exp(C_7\,t)\right) \leq \exp(B\,t).$$
Then, for $t_0$ small enough,
$$X(t_0) \leq C_8,$$
with $C_8$ depending on $t_0$, $B$ and $C_7$.

Now the theorem is almost proved for $f$ bounded. We have
\begin{multline*}
\mu\{x\in Q:|f(x)-f_Q|>\lambda\,\|f\|_*/t_0\} \\
\begin{split}
& \leq \int_Q
\exp\left(\frac{t_0\,|f(x)-f_Q|}{\|f\|_*}\right)\,\exp(-\lambda)\,d\mu(x)\\
& \leq C_8\,\mu(2 Q)\,\exp(-\lambda),
\end{split}
\end{multline*}
which is equivalent to \rf{mn88}.

When $f$ is not bounded, we consider the function $f_q$ of Lemma \ref{lemaq}.
By this Lemma and the subsequent remark we know that
$$\mu\{x\in Q:|f_q(x)-f_{Q,q}|>\lambda\} \leq C_5\,\mu(\rho
Q)\,\exp\left({\frac{-C_6\,\lambda}{\|f\|_*}}\right).$$
Since $\mu\{x\in Q:|f_q(x)-f_{Q,q}|>\lambda\} \to 
\mu\{x\in Q:|f(x)-f_{Q}|>\lambda\}$ as $q\to\infty$, 
\rf{mn88} holds in this case too.
\end{proof}

 From Theorem \ref{JN} we can get easily that the following spaces 
$\rbmo^p(\mu)$
coincide for all $p\in[1,\infty)$.
Given $\rho>1$ and $p\in[1,\infty)$, $\rbmo^p(\mu)$ is defined as follows.
We say that $f\in L^1_{loc}(\mu)$ is in
$\rbmo^p(\mu)$ if there exists some constant
$C_9$ such that for any cube $Q$ (centered at some point of $\supp(\mu)$)
\begin{equation} \label{def1p}
\left(\frac{1}{\mu(\rho Q)}\int_Q |f-m_{\wt{Q}}f|^p\,d\mu\right)^{1/p} \leq
C_9
\end{equation}
and
\begin{equation} \label{def2p}
|m_Q f-m_R f|\leq C_9 \,K_{Q,R} \quad \mbox{for any two
doubling cubes $Q\subset R$}.
\end{equation}
The minimal constant $C_9$ is the $\rbmo^p(\mu)$ norm of $f$, denoted by
$\|\cdot\|_{*,p}$.
Arguing as for $p=1$, one can show that another equivalent definition for
$\rbmo^p(\mu)$ can be given in terms of the numbers $\{f_Q\}_Q$
(as in \rf{def1'} and \rf{def2'}), and that the definition of the space does
not depend on the constant $\rho>1$.

We have the following corollary of John-Nirenberg inequality:

\begin{coro}
For $p\in[1,\infty)$, the spaces $\rbmo^p(\mu)$, coincide, and the norms
$\|\cdot\|_{*,p}$ are equivalent.
\end{coro}

\begin{proof}
The conditions \rf{def2} and \rf{def2p} coincide. So we only have
to compare \rf{def1} and \rf{def1p}.

For any $f\in L^1_{loc}(\mu)$, the inequality
$\|f\|_{*} \leq \|f\|_{*,p}$ follows from H\"older's inequality.
To obtain the converse inequality we will apply John-Nirenberg. If
$f\in\rbmo(\mu)$, then
\begin{eqnarray*}
\frac{1}{\mu(\rho Q)}\int_Q |f-m_{\wt{Q}}f|^p\,d\mu & = &
\frac{1}{\mu(\rho Q)}\int_0^\infty
p\,\lambda^{p-1}\,\mu\{x:|f(x)-m_{\wt{Q}}f|>\lambda\}\,d\lambda \\
& \leq &
C_5\, p\,\int_0^\infty \lambda^{p-1}\,
\exp\left({\frac{-C_6\,\lambda}{\|f\|_*}}\right)\,d\lambda \leq C\,\|f\|_*^p,
\end{eqnarray*}
and so
$\|f\|_{*,p} \leq C \|f\|_{*}$.
\end{proof}


\section{The space $\hba(\mu)$}

For a fixed $\rho>1$, a function $b\in L^1_{loc}(\mu)$ is called an {\em
atomic block} if
\begin{enumerate}
\item[1.] there exists some cube $R$ such that $\supp(b)\subset R$,
\item[2.] $\ds \int b\, d\mu = 0$,
\item[3.] there are functions $a_j$ supported on cubes $Q_j\subset R$ and
numbers $\lambda_j\in\R$ such that
$b=\sum_{j=1}^\infty \lambda_j a_j,$
and
$$\|a_j\|_{L^\infty(\mu)} \leq \left(\mu(\rho Q_j)\,K_{Q_j,R}\right)^{-1}.$$
\end{enumerate}
Then we denote
$$|b|_\hbm = \sum_j |\lambda_j|$$
(to be rigorous, we should think that $b$ is not only a function, but a
structure
formed by the function $b$, the cubes $R$ and $Q_j$, the functions $a_j$,
etc.)

Then, we say that $f\in \hbm$ if there are atomic blocks $b_i$
such that
$$f= \sum_{i=1}^\infty b_i,$$
with $\sum_i |b_i|_\hbm<\infty$. The $\hbm$ norm of $f$ is
$$\|f\|_\hbm = \inf \sum_i |b_i|_\hbm,$$
where the infimum is taken over all the possible decompositions of $f$
in atomic blocks.

Observe the difference with the atomic space $H^{1,\infty}_{at,\rho}(\mu)$.
The size condition on the functions $a_j$ is similar (we should forget the
coefficient $K_{Q,R}$), but the cancellation property $\int a_j\,d\mu=0$ is
substituted by something which offers more possibilities: We can gather
some terms $\lambda_j a_j$ in an atomic block $b$, and then we must have
$\int b\, d\mu=0$.

Notice also that if we take atomic blocks $b_i$ made up of a unique function
$a_i$, we derive $H^{1,\infty}_{at,\rho}(\mu) \subset \hbm$.

We have the following properties:

\begin{propo}  \label{prope2}
\begin{enumerate}
\item[1.] $\hbm$ is a Banach space.
\item[2.] $\hbm\subset L^1(\mu)$, with $\|f\|_{L^1(\mu)} \leq \|f\|_\hbm$.
\item[3.] The definition of $\hbm$ does not depend on the constant $\rho>1$.
\end{enumerate}
\end{propo}

\begin{proof}
The proofs of properties 1 and 2 are similar to the usual proofs for
$H^{1,\infty}_{at}(\mu)$.

Let us sketch the
proof of the third property, we can follow an argument similar
to the one of Lemma \ref{normeq1}. Given $\rho>\eta>1$, it is obvious that
$H^{1,\infty}_{atb,\rho}(\mu) \subset H^{1,\infty}_{atb,\eta}(\mu)$
and $\|f\|_{H^{1,\infty}_{atb,\eta}(\mu)} \leq
\|f\|_{H^{1,\infty}_{atb,\rho}(\mu)}$. For the converse inequality, given
an atomic block $b=\sum_j \lambda_j a_j$ with $\supp(a_j)\subset Q_j\subset
R$, it is not difficult to see that each function $a_j$ can be decomposed in a
finite fixed number of functions $a_{j,k}$ such that
$\|a_{j,k}\|_{L^\infty(\mu)} \leq \|a_{j}\|_{L^\infty(\mu)}$ for all $k$,
with $\supp(a_{j,k})\subset Q_{j,k}$, where $Q_{j,k}$ are cubes such that
$l(Q_{j,k}) \approx l(Q_j)$ and
$\rho Q_{j,k} \subset \eta Q_j$, etc.

Then, we will have $|b|_{H^{1,\infty}_{atb,\rho}(\mu)} \leq C\,
|b|_{H^{1,\infty}_{atb,\eta}(\mu)}$, which yields
$\|f\|_{H^{1,\infty}_{atb,\rho}(\mu)} \leq C\,
\|f\|_{H^{1,\infty}_{atb,\eta}(\mu)}. $
\end{proof}

Unless otherwise stated, we will assume that the constant $\rho$
in the definition $\hbm$ is equal to $2$.

Now we are going to see that if a CZO is bounded on $L^2(\mu)$, then
it is bounded from $\hbm$ into $L^1(\mu)$. In fact, we will
replace the assumption of $L^2(\mu)$ boundedness by another weaker assumption
(as in Theorem \ref{aco1}).

\begin{theorem} \label{aco2}
If for any cube $Q$ and any function $a$ supported on $Q$
\begin{equation}  \label{hip2}
\int_Q |T_\ve a|\,d\mu \leq C\,\|a\|_{L^\infty(\mu)}\,\mu(\rho Q)
\end{equation}
uniformly on $\ve>0$, then $T$ is bounded from $\hbm$ into $L^1(\mu)$.
\end{theorem}

\begin{proof} By standard arguments, it is enough to show that
$$\|T_\ve b\|_{L^1(\mu)} \leq C\,|b|_\hbm$$
for any atomic block $b$
with $\supp(b) \subset R$, $b=\sum_j \lambda_j\,a_j$, where the $a_j$'s are
functions satisfying the properties 3 and 4 of the definition of atomic
block.

We write
\begin{equation}  \label{hhj1}
\int |T_\ve b|\,d\mu = \int_{\R^d\setminus 2R} |T_\ve b|\,d\mu +
\int_{2R} |T_\ve b|\,d\mu.
\end{equation}
To estimate the first integral on the right hand side, we take into account
that $\int b\,d\mu=0$, and by usual arguments we get
\begin{equation}  \label{hhj3}
\int_{\R^d\setminus 2R} |T_\ve b|\,d\mu \leq C\,\|b\|_{L^1(\mu)}
\leq C\, |b|_\hbm.
\end{equation}

On the other hand, for the last integral in \rf{hhj1}, we have
\begin{eqnarray*}
\int_{2R} |T_\ve b|\,d\mu & \leq & \sum_j |\lambda_j| \int_{2R} |T_\ve
a_j|\,d\mu \\
& = & \sum_j |\lambda_j|
\int_{2Q_j} |T_\ve a_j|\,d\mu  +  \sum_j |\lambda_j| \int_{2R\setminus 2Q_j}
|T_\ve a_j|\,d\mu.
\end{eqnarray*}
By \rf{hip2}, for each $j$ we have
$$\int_{2Q_j} |T_\ve a_j|\,d\mu \leq C\,\|a_j\|_{L^\infty(\mu)}\,\mu(2\rho
Q_j).$$
Also,
\begin{eqnarray*}
\int_{2R\setminus 2Q_j} |T_\ve a_j|\, d\mu & \leq &
\sum_{k=1}^{N_{Q_j,R}} \int_{2^{k+1}Q_j\setminus 2^k Q_j}|T_\ve a_j|\, d\mu
\\
& \leq & C\sum_{k=1}^{N_{Q_j,R}}
\frac{\mu(2^{k+1}Q_j)}{l(2^{k+1}Q_j)^n}\,\|a_j\|_{L^1(\mu)} \\
& \leq & C\,K_{Q_j,R}\, \|a_j\|_{L^\infty(\mu)}\,\mu(Q_j).
\end{eqnarray*}
Thus
$$\int_{2R} |T_\ve a_j|\,d\mu \leq C\,
K_{Q_j,R}\, \|a_j\|_{L^\infty(\mu)}\,\mu(2\rho Q_j),$$
and then, taking into account the property 3 in Proposition \ref{prope2},
\begin{equation} \label{hhj4}
\int_{2R} |T_\ve b|\,d\mu \leq C\, \sum_j |\lambda_j|\,
K_{Q_j,R}\, \|a_j\|_{L^\infty(\mu)}\,\mu(2\rho Q_j) \leq C\, |b|_\hbm.
\end{equation}
By \rf{hhj3} and \rf{hhj4}, we are done. \end{proof}

We will see below, in Theorem \ref{equi}, that condition \rf{hip2} holds if
and only if $T$ is bounded from $\hbm$ into $L^1(\mu)$.

The spaces $\hbm$ and $\rbmo(\mu)$ are closely related. In the next section
we will prove that the dual of $\hbm$ is $\rbmo(\mu)$. Let us see one of the
inclusions (the easiest one).

\begin{lemma} \label{inclu1}
$$\rbmo(\mu) \subset \hbm^{\!\ast}.$$
That is, for $g\in \rbmo(\mu)$, the linear functional
$$L_g(f) = \int f\,g\,d\mu$$
defined over bounded functions $f$ with compact support extends to a
continuous linear functional $L_g$ over $\hbm$ with
$$\|L_g\|_{\hbm^{\!\ast}} \leq C\,\|g\|_*.$$
\end{lemma}

\begin{proof} Following some standard arguments (see \cite[p.294-296]{GR}, for
example), we only need to show that if  $b$ is an atomic block
and $g\in \rbmo(\mu)$, then
$$\left|\int b \,g\,d\mu\right| \leq C\,|b|_{\hbm}\,\|g\|_*.$$
Suppose $\supp(b) \subset R$, $b=\sum_j \lambda_j\,a_j$, where the $a_j$'s are
functions satisfying the properties 3 and 4 of the definition of atomic
blocks. Then, using $\int b\,d\mu=0$,
\begin{equation} \label{eqatb1}
\left|\int b \,g\,d\mu\right| =
\left|\int_R b \,(g-g_R)\,d\mu\right| \, \leq\, \sum_j
|\lambda_j|\, \|a_j\|_{L^\infty(\mu)}\, \int_{Q_j} |g-g_R|\,d\mu.
\end{equation}
Since $g\in \rbmo(\mu)$, we have
\begin{eqnarray*}
\int_{Q_j} |g-g_{R}|\,d\mu & \leq &
\int_{Q_j} |g-g_{Q_j}|\,d\mu + |g_R- g_{Q_j}|\,\mu(Q_j) \\
& \leq & \|g\|_*\,\mu(2Q_j) + K_{Q_j,R}\,\|g\|_*\,\mu(Q_j) \\
& \leq & C\,K_{Q_j,R}\,\|g\|_*\,\mu(2Q_j).
\end{eqnarray*}
 From \rf{eqatb1} we get
$$\left|\int b \,g\,d\mu\right| \leq C\,\sum_j |\lambda_j|\, \|g\|_* =
C\,|b|_{\hbm}\,\|g\|_*.$$
\end{proof}

In the following lemma we prove the converse inequality to the one stated
in Lemma \ref{inclu1}.

\begin{lemma} \label{inclu11}
If $g\in \rbmo(\mu)$, we have
$$\|L_g\|_{\hbm^{\!\ast}} \approx \|g\|_*.$$
\end{lemma}

To obtain this result we need to define another equivalent norm for
$\rbmo(\mu)$. First we introduce some notation.
 Given a cube $Q\subset \R^d$ and $f\in L^1_{loc}(\mu)$, let
$\alpha_Q(f)$ be the constant for which $\inf_{\alpha\in \R}
m_Q(|f-\alpha|)$ is attained. It is known that the constant $\alpha_Q(f)$
(which may be not unique) satisfies
$$\mu\{x\in Q:\,f(x)>\alpha_Q(f) \} \leq \frac{1}{2}\mu(Q)$$
and
$$\mu\{x\in Q:\,f(x)<\alpha_Q(f) \} \leq \frac{1}{2}\mu(Q)$$
(see \cite[p.30]{Journe}, for example).

Given $f\in L^1_{loc}(\mu)$, we denote by $\|f\|_\circ$ the minimal constant
$C_{10}$ such that
\begin{equation} \label{qqq2}
\frac{1}{\mu(2Q)}\, \int_Q |f-\alpha_\wt{Q}(f)|\,d\mu \leq C_{10}
\end{equation}
and, for any two doubling cubes $Q\subset R$,
\begin{equation} \label{qqq3}
|\alpha_{Q}(f) - \alpha_{R}(f)| \leq C_{10} \,K_{Q,R}.
\end{equation}
Then,

\begin{lemma} \label{normeq3}
$\|\cdot \|_\circ$ is a norm which equivalent with $\|\cdot\|_*$.
\end{lemma}

To prove this result, one can argue as in Lemma \ref{normeq2} and show that
the norm $\|\cdot\|_\circ$ is equivalent with the norm $\|\cdot\|_{**}$.
The details are left for the reader

\vspace{5mm}
\begin{proof}[Proof of Lemma \ref{inclu11}]
We have to prove that
$\|L_g\|_{\hbm^{\!\ast}} \geq C^{-1}\,\|g\|_*.$
We will show that there exists some function $f\in \hbm$ such that
$$|L_g(f)|\geq C^{-1}\, \|g\|_\circ \,\|f\|_\hbm.$$

Let $\ve>0$ be some small constant which will be fixed later. There are two
possibilities:
\begin{enumerate}
\item[1.] There exists some doubling cube $Q\subset \R^d$ such that
\begin{equation} \label{qqq0}
\int_Q |g-\alpha_Q(g)|\,d\mu \geq \ve\, \|g\|_\circ\,\mu(Q).
\end{equation}
\item[2.] For any doubling cube $Q\subset \R^d$, \rf{qqq0} does not hold.
\end{enumerate}

If case 1 holds and $Q$ is doubling and satisfies \rf{qqq0},
then we take $f$ such that $f(x)=1$ if $g(x)>\alpha_Q(g)$,
$f(x) = -1$ if $g(x)<\alpha_Q(g)$, and $f(x)=\pm 1$ if $g(x)=\alpha_Q(g)$, so
that $\int f\,d\mu=0$ (this is possible because of \rf{qqq2} and \rf{qqq3}).
Then,
$$\left| \int g\,f\,d\mu \right| =
\left| \int (g-\alpha_Q(g))\,f\,d\mu \right| =
\int |g-\alpha_Q(g)|\,d\mu \geq \ve\,\|g\|_\circ\,\mu(Q).$$
Since $f$ is an atomic block and $Q$ is doubling, $\|f\|_\hbm \leq
|f|_\hbm \leq C\,\mu(Q)$. Therefore
$$|L_g(f)| = \left| \int g\,f\,d\mu \right| \geq
C^{-1}\,\ve\,\|g\|_\circ\,\|f\|_\hbm.$$

In the case 2 there are again two possibilities:
\begin{enumerate}
\item[a)] For any two doubling cubes $Q\subset R$.
$$|\alpha_{Q}(g) - \alpha_{R}(g)|\leq \frac{1}{2}\,K_{Q,R}\,\|g\|_\circ.$$
\item[b)] There are doubling cubes $Q\subset R$ such that
$$|\alpha_{Q}(g) - \alpha_{R}(g)| > \frac{1}{2}\,K_{Q,R}\,\|g\|_\circ.$$
\end{enumerate}

Assume first that a) holds. By the definition of $\|g\|_\circ$ there
exists some cube $Q$ such that
$$\int_Q |g-\alpha_\wt{Q}(g)|\, d\mu \geq \frac{1}{2}\,\|g\|_\circ\,\mu(2Q).$$
We consider the following atomic block supported on $\wt{Q}$: We set
$f=a_1+a_2$, where
$$a_1 = \chi_{Q\cap\{g>\alpha_\wt{Q}(g)\} } -
\chi_{Q\cap\{g\leq \alpha_\wt{Q}(g)\} },$$
and $a_2$ is supported on $\wt{Q}$, constant on this cube, and such that
$\int (a_1+a_2)\,d\mu =0$.

Let us estimate $\|f\|_\hbm$. We have
$$\|a_2\|_{L^\infty(\mu)} \,\mu(\wt{Q}) = \left|\int a_2\,d\mu\right| =
\left|\int a_1\,d\mu\right| \leq \mu(Q).$$
Then, since $\wt{Q}$ is doubling and $K_{Q,2\wt{Q}}\leq C$,
$$\|f\|_\hbm \leq \|a_1\|_{L^\infty(\mu)}\, \mu(2Q) +
C\,\|a_2\|_{L^\infty(\mu)}\, \mu(2\wt{Q}) \leq C\,\mu(2Q).$$

Now we have
\begin{eqnarray}  \label{qqq21}
L_g(f) & = & \int g\,f\,d\mu = \int_\wt{Q} (g-\alpha_{\wt{Q}}(g)) \,f\,d\mu
\nonumber \\
& = & \int_\wt{Q} (g-\alpha_{\wt{Q}}(g)) \,a_1\,d\mu +
\int_\wt{Q} (g-\alpha_{\wt{Q}}(g)) \,a_2\,d\mu.
\end{eqnarray}
By the definition of $a_1$,
\begin{equation}  \label{qqq22}
\left|\int_\wt{Q} (g-\alpha_{\wt{Q}}(g)) \,a_1\,d\mu \right|
= \int_Q |g-\alpha_{\wt{Q}}(g)| \,d\mu
\geq \frac{1}{2}\,\|g\|_\circ\,\mu(2Q).
\end{equation}
On the other hand, by the computation about $\|a_2\|_{L^\infty(\mu)}$ and
since \rf{qqq0} does not hold for $\wt{Q}$,
\begin{equation}  \label{qqq23}
\left|\int_\wt{Q} (g-\alpha_{\wt{Q}}(g)) \,a_2\,d\mu\right| \leq
\frac{\mu(Q)}{\mu(\wt{Q})} \int_\wt{Q} |g-\alpha_{\wt{Q}}(g)| \,d\mu
\leq C\,\ve \|g\|_\circ\,\mu(2Q).
\end{equation}
By \rf{qqq21}, \rf{qqq22} and \rf{qqq23}, if $\ve$ has been chosen small
enough,
$$|L_g(f)|\geq \frac{1}{4}\,\|g\|_\circ\,\mu(2Q) \geq C^{-1}\,\|g\|_\circ\,
\|f\|_\hbm.$$

Now we consider the case b). Let $Q\subset R$ be doubling cubes such that
\begin{equation}  \label{qqq30}
|\alpha_{Q}(g) - \alpha_{R}(g)| > \frac{1}{2}\,K_{Q,R}\,\|g\|_\circ.
\end{equation}
We take
$$f = \frac{1}{\mu(R)}\,\chi_R - \frac{1}{\mu(Q)}\,\chi_Q.$$
So $\int f\,d\mu = 0$, and $f$ is an atomic block. Since $Q$ and $R$ are
doubling,
$\|f\|_\hbm \leq C\,K_{Q,R}.$
We have
\begin{eqnarray*}
L_g(f) & = & \int_R (g-\alpha_R(g))\,f\,d\mu \\
& = & \frac{1}{\mu(R)}\,
\int_R (g-\alpha_R(g))\,d\mu - \frac{1}{\mu(Q)}
\int_Q (g-\alpha_R(g))\,d\mu \\
& = & \frac{1}{\mu(R)}\, \int_R (g-\alpha_R(g))\,d\mu -
\frac{1}{\mu(Q)} \int_Q (g-\alpha_Q(g))\,d\mu \\
&& \mbox{}+ (\alpha_Q(g)- \alpha_R(g)).
\end{eqnarray*}
Since we are in the case 2, the terms
$$\left|\frac{1}{\mu(R)} \int_R (g-\alpha_R(g))\,d\mu\right|,\qquad
\left|\frac{1}{\mu(Q)} \int_Q (g-\alpha_Q(g))\,d\mu\right|$$
are bounded by $\ve\, \|g\|_\circ$. By \rf{qqq30}, if $\ve$ is chosen $\leq
1/8$, then
$$|L_g(f)| \geq \frac{1}{4} K_{Q,R}\,\|g\|_\circ \geq C^{-1}\,\|g\|_\circ\,
\|f\|_\hbm.$$
\end{proof}


\section{The spaces $\hbap(\mu)$ and duality}

To study the duality between $\hbm$ and $\rbmo(\mu)$ we will follow the
scheme of \cite[p.34-40]{Journe}. We will introduce the atomic spaces
$\hbap(\mu)$, and we will prove that they coincide with $\hbm$ and that the
dual of $\hbm$ is $\rbmo(\mu)$ simultaneously.

For a fixed $\rho>1$ and $p\in(1,\infty)$, a function $b\in L^1_{loc}(\mu)$ is
called a {\em $p$-atomic block} if
\begin{enumerate}
\item[1.] there exists some cube $R$ such that $\supp(b)\subset R$,
\item[2.] $\ds \int b\, d\mu = 0$,
\item[3.] there are functions $a_j$ supported in cubes $Q_j\subset R$ and
numbers $\lambda_j\in\R$ such that
$b=\sum_{j=1}^\infty \lambda_j a_j,$
and
$$\|a_j\|_{L^p(\mu)} \leq \mu(\rho Q_j)^{1/p-1}\, K_{Q_j,R}^{-1}.$$
\end{enumerate}
We denote
$$|b|_\hbp = \sum_j |\lambda_j|$$
(as in the case of $\hbm$, to be rigorous we should think that $b$ is not only
a function, but a structure
formed by the function $b$, the cubes $R$ and $Q_j$, the functions $a_j$,
etc.)

Then, we say that $f\in \hbp$ if there are $p$-atomic blocks $b_i$
such that
$$f= \sum_{i=1}^\infty b_i,$$
with $\sum_i |b_i|_\hbp<\infty$. The $\hbp$ norm of $f$ is
$$\|f\|_\hbp = \inf \sum_i |b_i|_\hbp,$$
where the infimum is taken over all the possible decompositions of $f$
in $p$-atomic blocks.

We have properties analogous to the ones for $\hbm$:

\begin{propo}  \label{prope3}
\begin{enumerate}
\item[1.] $\hbp$ is a Banach space.
\item[2.] $\hbp\subset L^1(\mu)$, with $\|f\|_{L^1(\mu)} \leq \|f\|_\hbp$.
\item[3.] For $1< p_1\leq p_2\leq \infty$, we have the continuous
inclusion $H^{1,p_2}_{atb}(\mu) \subset H^{1,p_1}_{atb}(\mu)$.
\item[4.] The definition of $\hbp$ does not depend on the constant $\rho>1$.
\end{enumerate}
\end{propo}

The proof of these properties is similar to the proof of the properties in
Proposition \ref{prope2}.

As in the case of $\rbmo(\mu)$ and $\hbm$, we will assume that the constant
$\rho$ in the definition of the $\hbp$ is $\rho=2$.

The proof about the duality between $\hbm$ and $\rbmo(\mu)$ and the
coincidence of the spaces $\hbp$ has been split in several lemmas. The
first one is the following.

\begin{lemma} \label{inclup}
For $1<p<\infty$,
$$\rbmo(\mu) \subset \hbp^{\!\ast}.$$
That is, for $g\in \rbmo(\mu)$, the linear functional
$$L_g(f) = \int f\,g\,d\mu$$
defined over bounded functions $f$ with compact support extends to a unique
continuous linear functional $L_g$ over $\hbp$ with
$$\|L_g\|_{\hbp^{\!\ast}} \leq C\,\|g\|_*.$$
\end{lemma}

\begin{proof} This lemma is very similar to Lemma \ref{inclu1}.
We only need to show that if $b$ is a $p$-atomic block
and $g\in \rbmo(\mu)$, then
$$\left|\int b \,g\,d\mu\right| \leq C\,|b|_{\hbp}\,\|g\|_*.$$
Suppose $\supp(b) \subset R$, $b=\sum_j \lambda_j\,a_j$, where the $a_j$'s are
functions satisfying the properties 3 and 4 of the definition of $p$-atomic
block. Since $\int b\,d\mu=0$,
\begin{equation} \label{eqatb1p}
\left|\int b \,g\,d\mu\right| =
\left|\int_R b \,(g-g_R)\,d\mu\right| \, \leq\, \sum_j
|\lambda_j|\, \|a_j\|_{L^p(\mu)}\, \left(\int_{Q_j}
|g-g_R|^{p'}\,d\mu\right)^{1/p'}\!\!,
\end{equation}
where $p'= p/(p-1)$.
As $g\in \rbmo(\mu) = \rbmo^{p'}(\mu)$, we have
\begin{eqnarray*}
\left(\int_{Q_j} |g-g_{R}|^{p'}\,d\mu\right)^{1/p'} & \leq &
\left(\int_{Q_j} |g-g_{Q_j}|\,d\mu\right)^{1/p'} + |g_R- g_{Q_j}|\,
\mu(Q_j)^{1/p'} \\
& \leq & \|g\|_*\,\mu(2Q_j)^{1/p'} + K_{Q_j,R}\,\|g\|_*\,\mu(Q_j)^{1/p'} \\
& \leq & C\,K_{Q_j,R}\,\|g\|_*\,\mu(2Q_j)^{1/p'}.
\end{eqnarray*}
 From \rf{eqatb1p} we get
$$\left|\int b \,g\,d\mu\right| \leq C\,\sum_j |\lambda_j|\, \|g\|_* =
C\,|b|_{\hbp}\,\|g\|_*.$$
\end{proof}

\begin{lemma} \label{pasclau}
For $1<p<\infty$,
$$\hbp^{\!\ast} \cap L^{p'}_{loc}(\mu) \subset \rbmo(\mu).$$
\end{lemma}

\begin{proof} Let $g \in L^{p'}_{loc}(\mu)$ such that such that the functional
$L_g$ belongs to $\hbp^{\!\ast}$. We have to show that $g\in \rbmo(\mu)$ and
$\|g\|_\circ \leq C\, \|L_g\|_{\hbp^{\!\ast}}$.
So we will see that, for any cube $Q$,
\begin{equation}  \label{qw1}
\frac{1}{\mu(2Q)} \int_Q |g-\alpha_\wt{Q}(g)|\,d\mu \leq C\, \|L_g\|_{\hbp^{\!\ast}}
\end{equation}
and for any two doubling cubes $Q\subset R$,
\begin{equation}  \label{qw2}
|\alpha_Q(g) - \alpha_R(g)|\leq C\,\|L_g\|_{\hbp^{\!\ast}}\,K_{Q,R}.
\end{equation}

First we will show that \rf{qw1} holds for any {\em doubling} cube $Q=\wt{Q}$.
In this case the argument is almost the same as the one of
\cite[p.38-39]{Journe}. We will repeat it for the sake of completeness.
Without loss of generality we may assume that
$$\int_{Q\cap\{g>\alpha_Q(g)\}} |g-\alpha_Q(g)|^{p'}\,d\mu \geq
\int_{Q\cap\{g<\alpha_Q(g)\}} |g-\alpha_Q(g)|^{p'}\,d\mu.$$
We consider an atomic block defined as follows:
$$
a(x) = \left\{ \begin{array}{ll}
|g(x) - \alpha_Q|^{p'-1}\quad & \mbox{if $x\in Q\cap\{g>\alpha_Q(g)\}$,} \\
C_{Q} & \mbox{if $x\in Q\cap\{g\leq\alpha_Q(g)\}$,} \\
0 & \mbox{if $x\not\in Q$,}
\end{array} \right.
$$
where $C_{Q}$ is a constant such that $\int a\,d\mu=0$.

By the definition of $\alpha_Q(g)$, we have
$$\mu(Q\cap\{g>\alpha_Q(g)\}) \leq \frac{1}{2}\,\mu(Q) \leq
\mu(Q\cap\{g\leq\alpha_Q(g)\}).$$
Since $Q$ is doubling,
\begin{multline*}
\|a\|_\hbp  \leq  C\,\|a\|_{L^p(\mu)} \,\mu(Q)^{1-1/p}
\leq  C\,\mu(Q)\\ \times \left(\frac{1}{\mu(Q)}
\int_{Q\cap\{g>\alpha_Q(g)\}} |g-\alpha_Q(g)|^{p'}\,d\mu + \frac{1}{\mu(Q)}
\int_{Q\cap\{g\leq\alpha_Q(g)\}} |C_{Q}|^p\,d\mu \right)^{1/p}.
\end{multline*}
Now we have
\begin{multline*}
\frac{1}{\mu(Q)} \int_{Q\cap\{g\leq\alpha_Q(g)\}} |C_{Q}|^p\,d\mu \\
\begin{split}
& \leq
\frac{1}{\mu(Q\cap \{g\leq \alpha_Q(g)\})}
\int_{Q\cap\{g\leq\alpha_Q(g)\}} |C_{Q}|^p\,d\mu \\
& =
\left|\frac{1}{\mu(Q\cap \{g\leq \alpha_Q(g)\})}
\int_{Q\cap\{g\leq\alpha_Q(g)\}} C_Q\,d\mu\right|^p \\
& =
\left(\frac{1}{\mu(Q\cap \{g\leq \alpha_Q(g)\})}
\int_{Q\cap\{g>\alpha_Q(g)\}} |g-\alpha_Q(g)|^{p'-1}\,d\mu\right)^p \\
& \leq
\frac{1}{\mu(Q\cap \{g\leq \alpha_Q(g)\})}
\int_{Q\cap\{g>\alpha_Q(g)\}} |g-\alpha_Q(g)|^{p'}\,d\mu.
\end{split}
\end{multline*}
Therefore,
\begin{equation}  \label{qww2}
\|a\|_\hbp \leq C\,\mu(Q)
\left(\frac{1}{\mu(Q)}
\int_{Q\cap\{g>\alpha_Q(g)\}} |g-\alpha_Q(g)|^{p'}\,d\mu \right)^{1/p}.
\end{equation}
Since $(g-\alpha_Q(g))a \geq 0$ on $Q$, we have
\begin{eqnarray}  \label{qww3}
\int_Q g\,a\,d\mu & = &
\int_Q (g-\alpha_Q(g))a\,d\mu \nonumber \\
& \geq & \int_{Q\cap\{g>\alpha_Q(g)\}} |g-\alpha_Q(g)|^{p'}\,d\mu
\geq \frac{1}{2} \int_Q |g-\alpha_Q(g)|^{p'}\,d\mu.
\end{eqnarray}
 From \rf{qww2} and \rf{qww3} we get
\begin{multline*}
\left(\frac{1}{\mu(Q)} \int_Q |g-\alpha_Q(g)|^{p'}\,d\mu\right)^{1/p'}\,
\|a\|_\hbp  \leq  C\, \int_Q |g-\alpha_Q(g)|^{p'}\,d\mu \\
\leq  C\,\int_Q g\,a\,d\mu
=  C\, L_g(a) \leq C\,\|L_g\|_{\hbp^{\!\ast}} \,\|a\|_\hbp.
\end{multline*}
So \rf{qw1} holds in this case.

\vspace{4mm}
Assume now that $Q$ is {\em non doubling}. We consider an atomic block
$b=a_1+a_2$, with
\begin{equation}  \label{qq44}
a_1=\frac{|g-\alpha_\wt{Q}(g)|^{p'}}{g-\alpha_\wt{Q}(g)}\,\chi_{Q\cap
\{g\neq \alpha_\wt{Q}(g)\}}
\end{equation}
and
\begin{equation}  \label{qq55}
a_2 = C_\wt{Q}\,\chi_\wt{Q},
\end{equation}
where $C_{\wt{Q}}$ is such that $\int (a_1+a_2)\,d\mu=0$.

Let us estimate $\|b\|_\hbp$.
Since $\wt{Q}$ is doubling and $K_{Q,\wt{Q}}\leq C$,
\begin{equation}  \label{qqww21}
\|b\|_\hbp \leq C\,\left(\int_Q |g-\alpha_\wt{Q}(g)|^{p'}\,d\mu
\right)^{1/p}\, \mu(2Q)^{1/p'} + C\,|C_{\wt{Q}}|\,\mu(\wt{Q}).
\end{equation}
Since $\int b\,d\mu=0$, we have
\begin{eqnarray} \label{qqww22}
\mu(\wt{Q}) \,|C_\wt{Q}| & = & \left|\int a_1\,d\mu\right|  \leq
\int_Q |g-\alpha_\wt{Q}(g)|^{p'-1}\,d\mu\nonumber \\ & \leq &
\left(\int_Q |g-\alpha_\wt{Q}(g)|^{p'}\,d\mu\right)^{1/p}\, \mu(Q)^{1/p'}.
\end{eqnarray}
Thus
\begin{equation}  \label{qqww33}
\|b\|_\hbp \leq C\,\left(\int_Q |g-\alpha_\wt{Q}(g)|^{p'}\,d\mu
\right)^{1/p}\, \mu(2Q)^{1/p'}.
\end{equation}

As $\int b\,d\mu=0$, we also have
$$\int g\,b\,d\mu =
\int_\wt{Q} (g-\alpha_\wt{Q}(g))b\,d\mu =
\int_Q (g-\alpha_\wt{Q}(g))a_1\,d\mu + C_\wt{Q}
\int_\wt{Q} (g-\alpha_\wt{Q}(g))\,d\mu.$$
Therefore, taking into account that $\wt{Q}$ satisfies \rf{qw1}, and using
\rf{qqww22},
\begin{multline}  \label{qq555}
\int_Q |g-\alpha_\wt{Q}(g)|^{p'}\,d\mu =
\int_Q (g-\alpha_\wt{Q}(g))a_1\,d\mu \\
\begin{split}
& \leq  \left| \int g\,b\,d\mu \right| + |C_\wt{Q}|
\int_\wt{Q} |g-\alpha_\wt{Q}(g)|\,d\mu \\
& \leq  \|L_g\|_{\hbp^{\!\ast}}\, \|b\|_\hbp + C\,|C_\wt{Q}|
\,\|L_g\|_{\hbp^{\!\ast}}\,\mu(\wt{Q}) \\
& \leq  C\,\|L_g\|_{\hbp^{\!\ast}}\, \left[\|b\|_\hbp +
\left(\int_Q |g-\alpha_\wt{Q}(g)|^{p'}\,d\mu\right)^{1/p}\, \mu(Q)^{1/p'}
\right].
\end{split}
\end{multline}
By \rf{qqww33} we get
$$
\int_Q |g-\alpha_\wt{Q}(g)|^{p'}\,d\mu \leq
C\,\|L_g\|_{\hbp^{\!\ast}}\,
\left(\int_Q |g-\alpha_\wt{Q}(g)|^{p'}\,d\mu\right)^{1/p}\, \mu(2Q)^{1/p'}.
$$
That is,
$$\left(\frac{1}{\mu(2Q)} \int_Q
|g-\alpha_\wt{Q}(g)|^{p'}\,d\mu \right)^{1/p'}
\leq C\, \|L_g\|_{\hbp^{\!\ast}},$$
which implies \rf{qw1}.

\vspace{4mm}
Finally, we have to show that \rf{qw2} holds for doubling cubes $Q\subset R$.
We consider an atomic block $b=a_1+a_2$ similar to the one defined above.
We only change $\wt{Q}$ by $R$ in \rf{qq44} and \rf{qq55}:
$$a_1=\frac{|g-\alpha_R(g)|}{g-\alpha_R(g)}\,\chi_{Q\cap
\{g\neq \alpha_R(g)\}}$$
and
$$a_2 = C_R\,\chi_R,$$
where $C_R$ is such that $\int (a_1+a_2)\,d\mu=0$.
Arguing as in \rf{qqww21}, \rf{qqww22} and \rf{qqww33} (the difference is that
now
$Q$ and $R$ are doubling, and we do {\em not} have $K_{Q,R}\leq C$) we will
obtain
\begin{equation}  \label{wr3}
\|b\|_\hbp \leq
C\,K_{Q,R}\,\left(\int_Q |g-\alpha_\wt{Q}(g)|^{p'}\,d\mu
\right)^{1/p}\, \mu(2Q)^{1/p'}.
\end{equation}
As in \rf{qq555}, we get
\begin{multline*}
\int_Q |g-\alpha_R(g)|^{p'}\,d\mu
\leq \|L_g\|_{\hbp^{\!\ast}}\, \|b\|_\hbp + C\,|C_R|
\,\|L_g\|_{\hbp^{\!\ast}}\,\mu(R) \\
\leq C\,\|L_g\|_{\hbp^{\!\ast}}\, \left[\|b\|_\hbp +
\left(\int_Q |g-\alpha_R(g)|^{p'}\,d\mu\right)^{1/p}\, \mu(Q)^{1/p'}
\right].
\end{multline*}
By \rf{wr3} we have
$$\int_Q |g-\alpha_R(g)|^{p'}\,d\mu \leq
C\,\|L_g\|_{\hbp^{\!\ast}}\, K_{Q,R}\,\left(\int_Q |g-\alpha_R(g)|^{p'}\,d\mu
\right)^{1/p}\, \mu(Q)^{1/p'}.$$
Thus
\begin{equation}  \label{fv1}
\left(\frac{1}{\mu(Q)} \int_Q
|g-\alpha_R(g)|^{p'}\,d\mu \right)^{1/p'}
\leq C\, \|L_g\|_{\hbp^{\!\ast}}\, K_{Q,R}.
\end{equation}
Since $Q$ is doubling and satisfies \rf{qw1}, by \rf{fv1} we get
\begin{eqnarray*}
|\alpha_Q(g)-\alpha_R(g)| & = & \frac{1}{\mu(Q)} \int_Q
|\alpha_Q(g)-\alpha_R(g)|\,d\mu \\
& \leq & \frac{1}{\mu(Q)} \int_Q  |g-\alpha_Q(g)|\,d\mu
+ \frac{1}{\mu(Q)} \int_Q  |g-\alpha_R(g)|\,d\mu \\
& \leq & C\, \|L_g\|_{\hbp^{\!\ast}}\, K_{Q,R},
\end{eqnarray*}
and we are done. \end{proof}

\begin{lemma} \label{pas2} For $1<p<\infty$,
$$\hbp^{\!\ast} \subset L^{p'}_{loc}(\mu) \hspace{3mm}\mbox{modulo constants.}$$
\end{lemma}

\begin{proof}
This lemma is an easy consequence of the Riesz representation theorem.
The same argument as the one of \cite[p.39-40]{Journe} works here.
\end{proof}

Now we have:

\begin{theorem}  \label{dualitat}
For $1<p<\infty$, $\hbp = \hbm$. Also, $\hbm^{\!\ast} = \rbmo(\mu)$.
\end{theorem}

\begin{proof}
Notice that, by Lemmas \ref{inclup}, \ref{pasclau} and \ref{pas2},
$\hbp^{\!\ast} = \rbmo(\mu)$ for $1<p<\infty$.

Now we repeat the arguments in \cite{Journe} again. We consider the diagram
$$
\begin{array}{rll}
i: \hbm & \lra & \hbp \\
i^{\ast}: \rbmo(\mu) = \hbp^{\!\ast} & \lra & \hbm^{\!\ast}.
\end{array}
$$
The map $i$ is an inclusion and $i^{\ast}$ is the canonical injection of
$\rbmo(\mu)$ in $\hbm^{\!\ast}$ (with the identification $g\equiv L_g$ for $g\in
\rbmo(\mu)$). By Lemma \ref{inclu11}, $i^{\ast}(\rbmo(\mu))$ is a closed subspace
of $\hbm^{\!\ast}$. By Banach's closed range theorem (see \cite[p.205]{Yosida}),
$\hbm$ is also closed in $\hbp$. Now it is easily checked that the Hahn-Banach
theorem and the fact that $\hbp^{\!\ast}=\rbmo(\mu)$ imply $\hbm^{\!\ast}=\rbmo(\mu)$.
\end{proof}

\bexam \label{ex4}
By the previous theorem and the fact that for an AD-regular set the space
$\rbmo(\mu)$ coincides with $\bmo(\mu)$, we derive that in this case we have
$\hbm= H^{1,\infty}(\mu)$ too (using the same sort of uniqueness argument as
above). However, this does not hold for all
doubling measures $\mu$. For instance, in the example \ref{ex2}
($\mu$ equal to the Lebesgue measure on $[0,1]^2$, with $n=1$, $d=2$)
since $\rbmo(\mu)=L^\infty(\mu)$ modulo constants, we have $\hbm=
\{f\in L^1(\mu):\,\int f\,d\mu=0\}$.
\eexam


\section{The sharp maximal operator}

The classical (centered) sharp maximal operator $M^\sharp$ is defined as
$$M^\sharp f(x) = \sup_{Q} \frac{1}{\mu(Q)} \int_Q |f-m_Qf|\,d\mu,$$
where the supremum is taken over the cubes $Q$ centered  at $x$ and
$f\in L^1_{loc}(\mu)$. Then, one has $f\in\bmo(\mu)$ if and only if
$M^\sharp f\in L^\infty(\mu)$.

Since $M^\sharp f$ is pointwise bounded above by the
(centered) Hardy-Littlewood maximal operator
$$Mf(x) = \sup_{Q} \frac{1}{\mu(Q)} \int_Q |f|\,d\mu$$
(with the supremum again over the cubes $Q$ centered at $x$),
one has
\begin{equation}  \label{sh1}
\|M^\sharp f\|_{L^p(\mu)} \leq C\,  \|f\|_{L^p(\mu)}
\end{equation}
for $1<p\leq \infty$.
On the other hand, the converse inequality also holds. If $f\in L^p(\mu)$,
$1<p<\infty$, then
\begin{equation}  \label{sh2}
\|f\|_{L^p(\mu)} \leq C\, \|M^\sharp f\|_{L^p(\mu)}
\end{equation}
(assuming $\int f\,d\mu=0$ if $\|\mu\|<\infty$).

In \cite{MMNO} it is shown that the inequalities \rf{sh1} and \rf{sh2} are
satisfied too if $\mu$ is non doubling (choosing an appropiate grid of cubes).
However, the above definition of the sharp operator
is not useful for our purposes because we do not have the equivalence
\begin{equation} \label{sh3}
f\in\rbmo(\mu) \iff M^\sharp f\in L^\infty(\mu).
\end{equation}

Now we want to introduce another sharp maximal operator suitable for our
space $\rbmo(\mu)$ enjoying properties similar to the ones of the classical
sharp operator. We define
\begin{equation} \label{defsharp}
M^{\sharp} f(x) = \sup_{Q\ni x} \frac{1}{\mu({\textstyle
\frac{3}{2}}Q)}
\int_Q |f- m_\wt{Q} f|\, d\mu + \!\!\!\sup_{
\mbox{\scriptsize $
\begin{array}{c} Q\subset R:\,x\in Q,\\ \mbox{$Q,R$ doubling}
\end{array}$}} \!\!\!
\frac{|m_Qf-m_R f|}{K_{Q,R}}.
\end{equation}
Notice that the cubes that appear in these supremums may be
non centered at $x$. It is clear that $f\in \rbmo(\mu)$ if and only if
$M^{\sharp} f\in L^\infty(\mu)$ (recall Remark \ref{remnou}).

We consider the non centered doubling maximal operator $N$:
$$N f(x) = \sup_{
\mbox{\scriptsize $
\begin{array}{c} Q\ni x,\\ \mbox{$Q$ doubling}
\end{array}$}}
\frac{1}{\mu(Q)} \int_Q |f|\, d\mu.$$
Observe that $|f(x)| \leq Nf(x)$ for $\mu$-a.e. $x\in \R^d$, by Remark
\ref{rem9}. Moreover,
the operator $N$ is of weak type $(1,1)$ and bounded on $L^p(\mu)$,
$p\in(1,\infty]$. Indeed, if $Q$ is doubling and $x\in Q$, we can
write
$$\frac{1}{\mu(Q)} \int_Q |f|\,d\mu \leq \frac{\beta_d}{\mu(2Q)}
\int_{Q}|f|\,d\mu \leq \beta_d\,M_{(2)}f(x),$$
where, for $\rho>1$, we denote
$$M_{(\rho)} f(x) = \sup_{Q\ni x} \frac{1}{\mu(\rho Q)}\int_Q |f|\,d\mu.$$
So $N f(x) \leq \beta_d M_{(2)}f(x)$.
The maximal operator $M_{(\rho)}$ is bounded above by the operator
defined as
$$M^{(\rho)} f(x) = \sup_{\rho^{-1}Q\ni x} \frac{1}{\mu(Q)}\int_Q
|f|\,d\mu.$$
This is the version of the Hardy-Littlewood operator that one
obtains taking supremums over cubes $Q$ which may be non centered at $x$
but such that $x\in \rho^{-1} Q$. Recall that since $0<\rho^{-1}<1$, one can
apply Besicovich's covering theorem (see \cite{Morse} or \cite[p.6-7]{Guzman}, 
for example) and then one gets that $M^{(\rho)}$ is of
weak type $(1,1)$ and bounded on $L^p(\mu)$, $p\in(1,\infty]$. As a
consequence, $M_{(\rho)}$ is also of
weak type $(1,1)$ and bounded on $L^p(\mu)$, $p\in(1,\infty]$

Now we derive that $M^{\sharp}$ also satifies the inequality \rf{sh1}
since the first supremum in the definition of $M^{\sharp} f$ is bounded
by $M_{(3/2)}f(x)+N f(x)$ while the second one is bounded by
$2N f(x)$.

\brem \label{rem40}
We have
$$M^{\sharp} |f|(x) \leq 5\beta_d\, M^{\sharp} f(x).$$
This easy to check: Assume that $x\in Q$ and $Q$ is
doubling. Then we have
\begin{eqnarray} \label{bb23}
\left|m_Q|f| - |m_Q f| \right| & = & \left|\frac{1}{\mu(Q)} \,\int_Q
(|f(x)|-|m_Q f|)\, d\mu(x) \right| \nonumber \\
& \leq & \frac{1}{\mu(Q)} \,\int_Q |f(x)-m_Q f| \, d\mu(x) \nonumber \\
& \leq & \beta_d \,M^{\sharp} f(x).
\end{eqnarray}
Therefore, if $Q\subset R$ are doubling,
\begin{eqnarray*}
\left|m_Q|f| - m_R|f|\, \right| & \leq & \left|m_Q f - m_R f \right|
+ 2\beta_d \,M^{\sharp} f(x) \\
& \leq & (K_{Q,R}+2\,\beta_d)\,M^{\sharp} f(x)
\leq 3\beta_d K_{Q,R}\,M^{\sharp} f(x).
\end{eqnarray*}
Thus
$$\sup_{
\mbox{\scriptsize $
\begin{array}{c} Q\subset R:\,x\in Q,\\ \mbox{$Q,R$ doubling}
\end{array}$}} \!\!\!
\frac{|m_Q|f|-m_R|f|\,|}{K_{Q,R}} \leq 3\beta_d \,M^{\sharp} f(x).$$
For the other supremum, by \rf{bb23} we have
\begin{multline*}
\sup_{Q\ni x} \frac{1}{\mu({\textstyle
\frac{3}{2}}Q)}
\int_Q |\,|f(x)|- m_\wt{Q} |f|\,|\, d\mu(x) \\
\begin{split}
& \leq
\sup_{Q\ni x} \frac{1}{\mu({\textstyle
\frac{3}{2}}Q)}
\int_Q |\,|f(x)|- |m_\wt{Q} f|\,|\, d\mu(x) + \beta_d \,M^{\sharp} f(x)\\
&\leq
\sup_{Q\ni x} \frac{1}{\mu({\textstyle
\frac{3}{2}}Q)}
\int_Q |f(x)- m_\wt{Q} f|\, d\mu(x) + \beta_d \,M^{\sharp} f(x) \\
& \leq  2\beta_d \,M^{\sharp} f(x).
\end{split}
\end{multline*}
\erem

Finally we are going to prove that our new operator $M^{\sharp}$
satisfies \rf{sh2} too. This is a consequence of
the next result and Remark \ref{rem9}.

\begin{theorem}  \label{sharp}
Let $f\in L^1_{loc}(\mu)$, with $\int f\,d\mu=0$ if $\|\mu\|<\infty$.
For $1<p<\infty$, if $\inf(1,N f)\in L^{p}(\mu)$,
then we have
\begin{equation}  \label{gli-1}
\|N f\|_{L^p(\mu)} \leq C\, \|M^{\sharp} f\|_{L^p(\mu)}.
\end{equation}
\end{theorem}

\begin{proof}
We assume $\|\mu\|=\infty$. The proof for $\|\mu\|<\infty$ is similar.
For some fixed $\eta<1$ and all $\ve>0$, we will prove that there exists some
$\delta>0$ such that for any $\lambda>0$ we have the following good
$\lambda$ inequality:
\begin{equation}  \label{gli0}
\mu\{x:\, N f(x)> (1+\ve)\lambda,\, M^{\sharp} f(x) \leq
\delta\lambda\} \leq \eta \,\mu\{x:\,N f(x)>\lambda\}.
\end{equation}
It is well known that by this inequality one gets
$\|N f\|_{L^p(\mu)} \leq C\, \|M^{\sharp} f\|_{L^p(\mu)}.$
if $\inf(1,N f)\in L^{p}(\mu)$.

We denote $\Omega_\lambda = \{x:\,N f(x)>\lambda\}$ and
$$E_\lambda = \{x:\, N f(x)> (1+\ve)\lambda,\, M^{\sharp} f(x)
\leq \delta\lambda\}.$$
For the moment we assume $f\in L^p(\mu)$. For each $x\in E_\lambda$,
among the doubling cubes $Q$
that contain $x$ and such that $m_Q|f|>(1+\ve/2) \lambda$,
we consider one cube $Q_x$ which has `almost maximal' side lenght, in the
sense that if some doubling
cube $Q'$ with side lenght $\geq 2l(Q_x)$ contains $x$, then
$m_{Q'}|f|\leq(1+\ve/2)\lambda$.
It is easy to check that this maximal cube $Q_x$ exists, because $f\in
L^p(\mu)$.

Let $R_x$ be the cube centered at $x$ with side length $3l(Q_x)$. We denote
$S_x= \wt{R_x}$. Then, assuming $\delta$ small enough we have
$m_{S_x}|f|>\lambda$, and then
$S_x\subset \Omega_\lambda$. Indeed, by construction, we have
$K_{Q_x,S_x}\leq C$. Then, as $Q_x\subset S_x$ are doubling cubes
containing $x$,
$$\left|m_{Q_x}|f| - m_{S_x}|f|\,\right|\leq K_{Q_x,S_x}\,
M^{\sharp} |f|(x) \leq C_{11}
5\beta_d\delta\lambda.$$
Thus, for $\delta<C_{11}\ve/10\beta_d$,
$$m_{S_x}|f|>(1+\ve/2)\lambda - C_{11}5\beta_d\delta\lambda >\lambda.$$

By Besicovich's covering theorem there are $n_B$ (depending on $d$)
subfamilies $\DD_k=\{S_i^k\}_i$, $k=1,\ldots,n_B$, of cubes $S_x$ such that
they cover $E_\lambda$, they are centered at points $x_i^k\in E_\lambda$,
and each subfamily $\DD_k$ is disjoint.
Therefore, at least one the subfamilies $\DD_k$ satisfies
$$\mu\left( \bigcup_i S_i^k \right) \geq \frac{1}{n_B}\, \mu\left(
\bigcup_{i,k} S_i^k \right).$$
Suppose, for example, that it is $\DD_1$.
We will prove that for each cube $S_i^1$,
\begin{equation}   \label{gli2}
\mu(S_i^1\cap E_\lambda) \leq \,\mu(S_i^1)/2n_B
\end{equation}
if $\delta$ is chosen small enough.
 From this inequality one gets
$$\mu\left(E_\lambda\cap\bigcup_{i} S_i^1\right) \leq
\frac{1}{2n_B}\,\sum_i \mu(S_i^1) \leq \frac{1}{2n_B}\, \mu(\Omega_\lambda).$$
Then,
\begin{eqnarray} \label{gli10}
\mu(E_\lambda) & \leq & \mu\left(\bigcup_{i,k} S_i^k \setminus
\bigcup_{i} S_i^1\right) + \mu\left(E_\lambda\cap\bigcup_{i} S_i^1\right)
\nonumber\\
& \leq & \left(1-\frac{1}{n_B}\right)\,
\mu\left(\bigcup_{i} S_i^1\right) + \frac{1}{2n_B}\,\mu(\Omega_\lambda)
\nonumber\\
& \leq & \left(1-\frac{1}{2n_B}\right) \,\mu(\Omega_\lambda).
\end{eqnarray}

Let us prove \rf{gli2}. Let $y\in S_i^1\cap E_\lambda$. If $Q\ni y$ is
doubling and such that $m_Q |f| > (1+\ve)\lambda$, then
$l(Q)\leq l(S^1_i)/8$. Otherwise, $\wt{30Q}\supset S^1_i\supset Q_{x_i^1}$,
and since $Q$ and $\wt{30Q}$ are doubling, we have
$$\left|m_Q|f| - m_\wt{30Q} |f|\right| \leq
K_{Q,\wt{30Q}}\,M^{\sharp}|f|(y) \leq C_{12}\,\delta\lambda
\leq \frac{\ve}{2}\,\lambda,$$
assuming $C_{12}\delta<\ve/2$, and so
$$m_\wt{30Q} |f| >(1+\ve/2)\,\lambda,$$
which contradicts the choice of $Q_{x_i^1}$ because $\wt{30Q}\supset
Q_{x_i^1}$ and $l(\wt{30Q})>2 l(Q_{x_i^1})$.

So $N f(y)>(1+\ve)\lambda$, implies
$$
N (\chi_{\frac{5}{4}S^1_i} f)(y)>(1+\ve)\lambda.
$$
On the other hand, we also have
$$m_{\wt{\frac{5}{4}S^1_i}} |f| \leq (1+\ve/2)\,\lambda,$$
since $\wt{\frac{5}{4}S^1_i}$ is doubling and its
side length is $>2l(Q_{x^1_i})$.
Therefore, we get
$$N (\chi_{\frac{5}{4}S^1_i} |f| - m_{\wt{\frac{5}{4}S^1_i}} |f|)
(y)>\frac{\ve}{2}\,\lambda,$$
and then, by the weak $(1,1)$ boundedness of $N$ and the fact that
$S^1_i$ is doubling,
\begin{eqnarray*}
\mu(S_i^1\cap E_\lambda) & \leq & \mu\{y:\,
N (\chi_{\frac{5}{4}S^1_i} (|f| - m_{\wt{\frac{5}{4}S^1_i}} |f|)
)(y)>\frac{\ve}{2}\,\lambda\} \\
& \leq & \frac{C}{\ve\lambda} \int_{\frac{5}{4}S^1_i}
(|f| - m_{\wt{\frac{5}{4}S^1_i}} |f|)\, d\mu \\
& \leq & \frac{C}{\ve\lambda} \, \mu(2S^1_i)\,M^\sharp |f|(x^1_i) \\
& \leq & \frac{C_{13}\,\delta}{\ve}\, \mu(S^1_i).
\end{eqnarray*}
Thus, \rf{gli2} follows by choosing $\delta<\ve/2n_B C_{13}$, which implies
\rf{gli10}, and as a consequence we obtain \rf{gli0} and \rf{gli-1} (under the
assumption $f\in L^p(\mu)$).

\vspace{2mm}
Suppose now that $f\not\in L^p(\mu)$. We consider the 
functions $f_q$, $q\geq1$, introduced in Lemma \ref{lemaq}. Since
for all functions $g,h\in L^1_{loc}(\mu)$ and all $x$ we have 
$M^\sharp(g+h)(x) \leq
M^\sharp g(x) + M^\sharp h(x)$ and $M^\sharp|g|(x) \leq C\,M^\sharp g(x)$,
operating as in Lemma \ref{lemaq} we get $M^\sharp f_q(x) \leq C\,M^\sharp
f(x)$. On the other hand, $|f_q(x)| \leq q\,\inf(1,|f|)(x) \leq 
q\,\inf(1,Nf)(x)$ and so $f_q\in L^p(\mu)$. Therefore,
$$\|N f_q\|_{L^p(\mu)} \leq C\, \|M^{\sharp} f\|_{L^p(\mu)}.$$
Taking the limit as $q\to\infty$, \rf{gli-1} follows.
\end{proof}


\section{Interpolation results}

An immediate corollary of the properties of the sharp operator is the
following result.

\begin{theorem} \label{interfacil}
Let $1<p<\infty$ and $T$ be a linear operator bounded on $L^p(\mu)$
and from $L^\infty(\mu)$ into $\rbmo(\mu)$. Then $T$ extends boundedly to
$L^r(\mu)$, $p<r<\infty$.
\end{theorem}

\begin{proof} Assume $\|\mu\|=\infty$.
The operator $M^{\sharp}\circ T$ is sublinear and it is bounded
in $L^p(\mu)$ and $L^\infty(\mu)$. By the Marcinkiewitz interpolation theorem,
it is bounded on $L^r(\mu)$, $p<r<\infty$. That is,
$$\|M^{\sharp}Tf\|_{L^r(\mu)} \leq C\,\|f\|_{L^r(\mu)}.$$
We may assume that $f\in L^r(\mu)$ has compact support. Then $f\in L^p(\mu)$
and so $Tf\in L^p(\mu)$. Thus
$N f\in L^p(\mu)$, and so $\inf(1,Nf)\in L^r(\mu)$.
By Theorem \ref{sharp}, we have
$$\|Tf\|_{L^r(\mu)}\leq C\,\|M^{\sharp}Tf\|_{L^r(\mu)}
\leq C\,\|f\|_{L^r(\mu)}.$$

The proof for $\|\mu\|<\infty$ is similar:
Given $f\in L^r(\mu)$, we write
$f= \left(f-\int f\,d\mu\right) + \int f\,d\mu$.
It easily seen that the same argument as for $\|\mu\|=\infty$
can be applied to the function $f-\int f\,d\mu$.
On the other hand,
$T$ is bounded on $L^r(\mu)$ over constant functions. Indeed,
since $T$ is bounded from $L^\infty(\mu)$ into $\rbmo(\mu)$, we get
\begin{eqnarray*}
\|T1\|_{L^r(\mu)} & \leq &
\|T1- m_{\R^d}(T1)\|_{L^r(\mu)} + \|m_{\R^d}(T1)\|_{L^r(\mu)} \\
& \leq & C\,\mu(\R^d)^{1/r} + \|m_{\R^d}(T1)\|_{L^r(\mu)}.
\end{eqnarray*}
We also have 
\begin{eqnarray*}
\|m_{\R^d}(T1)\|_{L^r(\mu)} & = & \|m_{\R^d}(T1)\|_{L^p(\mu)}\,
\mu(\R^d)^{\frac{1}{r}-\frac{1}{p}} \\
& \leq & \|T1\|_{L^p(\mu)}\, \mu(\R^d)^{\frac{1}{r}-\frac{1}{p}} \,\leq \,C\,
\mu(\R^d)^{1/r} \, \equiv \,C\,\|1\|_{L^r(\mu)},
\end{eqnarray*}
and so $\|T1\|_{L^r(\mu)} \leq C\,\|1\|_{L^r(\mu)}$.
\end{proof}

The main theorem of this section is another interpolation result which
is not as immediate as the previous one. Using this result, we will be able
to prove the $T(1)$ theorem for the Cauchy integral in the
next section. The statement is the following.

\begin{theorem} \label{int*}
Let $T$ be a linear operator which is bounded from $\hbm$ into $L^1(\mu)$
and from $L^\infty(\mu)$ into $\rbmo(\mu)$. Then, $T$ extends boundedly to
$L^p(\mu)$, $1<p<\infty$.
\end{theorem}

The proof of this theorem will follow the scheme of \cite[p.43-44]{Journe}.
To prove it we need a substitute for the Calder\'on-Zygmund
decomposition of a function, suitable for non doubling measures.
Nazarov, Treil and Volberg \cite{NTV2}
showed that if a CZO is bounded on
$L^2(\mu)$, then it is of weak type $(1,1)$ (with $\mu$ non doubling).
They used some kind of Calder\'on-Zygmund decomposition
to obtain this result. However, their decomposition does not work in the proof
of Theorem \ref{int*}. 
Mateu, Mattila, Nicolau and Orobitg \cite{MMNO} also used a Calder\'on-Zygmund
type decomposition to prove an interpolation theorem between 
$(H^1_{at}(\mu), L^1(\mu))$ and $(L^\infty(\mu),\bmo(\mu))$ 
with $\mu$ non doubling. Their decompositon is not suitable for our purposes
either. We will use the following decomposition instead.

\begin{lemma}[Calder\'on-Zygmund decomposition]
For $1\leq\! p< \!\infty$, consider $f\in L^p(\mu)$ with compact support.
For any $\lambda>0$ (with $\lambda>\beta_d\,\|f\|_{L^1(\mu)}/\|\mu\|$ if 
$\|\mu\|<\infty$), we have:
\begin{itemize}
\item[a)] There exists a finite family of
almost disjoint (i.e. with a bounded overlap) cubes $\{Q_i\}_i$ such that
\begin{equation}  \label{cc1}
\frac{1}{\mu(2Q_i)} \int_{Q_i} |f|^p\,d\mu >\frac{\lambda^p}{\beta_d},
\end{equation}
\begin{equation}  \label{cc2}
\frac{1}{\mu(2\eta Q_i)} \int_{\eta Q_i} |f|^p\,d\mu \leq
\frac{\lambda^p}{\beta_d} \quad
\mbox{for all $\eta >2$,}
\end{equation}
\begin{equation}  \label{cc3}
|f|\leq \lambda \quad \mbox{a.e. ($\mu$) on $\R^d\setminus\bigcup_i
Q_i$}.
\end{equation}

\item[b)] For each $i$, let $w_i= \frac{\chi_{Q_i}}{\sum_k \chi_{Q_k}}$
and let $R_i$ be a  $(6,6^{n+1})$-doubling cube concentric
with $Q_i$, with $l(R_i)>4l(Q_i)$. Then there exists a family of functions
$\vphi_i$ with $\supp(\vphi_i)\subset R_i$ satisfying
\begin{equation}  \label{cc4}
\int \vphi_i \,d\mu = \int_{Q_i} f\,w_i\,d\mu,
\end{equation}
\begin{equation}  \label{cc5}
\sum_i |\vphi_i| \leq B\,\lambda
\end{equation}
(where $B$ is some constant), and if $1<p<\infty$,
\begin{equation}  \label{cc6}
\left(\int_{R_i} |\vphi_i|^p \,d\mu\right)^{1/p}\mu(R_i)^{1/p'} \leq
\frac{C}{\lambda^{p-1}} \int_{Q_i}|f|^p\,d\mu.
\end{equation}

\item[c)] For $1<p<\infty$, if $R_i$ is the smallest
$(6,6^{n+1})$-doubling cube of
the family $\{6^k Q_i\}_{k\geq 1}$ and we set $b =
\sum_i (f\,w_i-\vphi_i)$, then
\begin{equation} \label{cc7}
\| b\|_\hbp \leq \frac{C}{\lambda^{p-1}} \|f\|_{L^p(\mu)}^p.
\end{equation}
\end{itemize}
\end{lemma}

\begin{proof}
We will assume $\|\mu\|=\infty$.
\begin{itemize}
\item[a) ] Taking into account Remark \ref{rem9},
for $\mu$-almost all $x\in\R^d$ such that $|f(x)|^p>\lambda^p$, there exists a
cube $Q_x$ satisfying
$$\frac{1}{\mu(2Q_x)} \int_{Q_x} |f|^p \,d\mu > \frac{\lambda^p}{\beta_d}$$
and such that if $Q_x'$ is centered at $x$ with $l(Q_x')>2l(Q_x)$, then
$$\frac{1}{\mu(2Q_x')} \int_{Q_x'} |f|^p \,d\mu \leq
\frac{\lambda^p}{\beta_d}.$$
Now we can apply Besicovich's covering theorem to get an almost disjoint 
subfamily of cubes
$\{Q_i\}_i \subset \{Q_x\}_x$ satisfying \rf{cc1}, \rf{cc2} and \rf{cc3}.

\item[b)] Assume first that the family of cubes $\{Q_i\}_i$ is finite.
Then we may suppose that this family of cubes is ordered in such a way
that the sizes of the cubes $R_i$ are non decreasing (i.e. $l(R_{i+1})\geq
l(R_i)$).
The functions $\vphi_i$ that we will construct will be of the form $\vphi_i
=\alpha_i\,\chi_{A_i}$, with $\alpha_i\in\R$ and $A_i\subset R_i$.
We set $A_1=R_1$ and
$$\vphi_1 = \alpha_1\,\chi_{R_1},$$
where the constant $\alpha_1$ is chosen so that $\int_{Q_1}f\,w_1\,d\mu=\int
\vphi_1\,d\mu$.

Suppose that $\vphi_1,\ldots,\vphi_{k-1}$ have been constructed,
satisfy \rf{cc4} and $\sum_{i=1}^{k-1} |\vphi_i|\leq
B\,\lambda,$  where $B$ is some constant (which will be fixed below).

Let $R_{s_1},\ldots,R_{s_m}$ be the subfamily of
$R_1,\ldots,R_{k-1}$ such that $R_{s_j}\cap R_k \neq \varnothing$.
As $l(R_{s_j}) \leq l(R_k)$ (because of the non decreasing sizes of $R_i$),
we have $R_{s_j} \subset 3R_k$. Taking into account that for $i=1,\ldots,k-1$
$$\int |\vphi_i|\,d\mu \leq \int_{Q_i} |f|\,d\mu$$
by \rf{cc4}, and using that $R_k$ is $(6,6^{n+1})$-doubling and \rf{cc2}, we
get
\begin{eqnarray*}
\sum_j \int |\vphi_{s_j}|\,d\mu & \leq & \sum_j
\int_{Q_{s_j}}|f|\,d\mu\\
& \leq & C \int_{3R_k} |f|\,d\mu \\
& \leq & C \left(\int_{3R_k} |f|^p\,d\mu\right)^{1/p}\,\mu(3R_k)^{1/p'}\\
& \leq & C \lambda\mu(6R_k)^{1/p}\,\mu(3R_k)^{1/p'}\\
& \leq & C_{14}\lambda\,\mu(R_k).
\end{eqnarray*}
Therefore,
$$\mu\left\{{\textstyle \sum_j} |\vphi_{s_j}| > 2C_{14}\lambda\right\}\leq
\frac{\mu(R_k)}{2}.$$ So we set
$$A_k = R_k\cap\left\{{\textstyle \sum_j} |\vphi_{s_j}| \leq
2C_{14}\lambda\right\},$$
and then $$\mu(A_k)\geq \mu(R_k)/2.$$

The constant $\alpha_k$ is chosen so that for $\vphi_k = \alpha_k\, \chi_{A_k}$
we have $\int\vphi_k\,d\mu = \int_{Q_k} f\,w_k\,
d\mu$.  Then we obtain
\begin{eqnarray*}  |\alpha_k|  &\leq&  \frac{1}{\mu(A_k)}\int_{Q_k}  |f|\,d\mu
\leq     \frac{2}{\mu(R_k)}\int_{Q_k}     |f|\,d\mu     \\     &     \leq    &
\frac{2}{\mu(R_k)}\int_{\frac{1}{2}R_k}             |f|\,d\mu             \leq
\left(\frac{2}{\mu(R_k)}\int_{\frac{1}{2}R_k} |f|^p\,d\mu\right)^{1/p}
 \leq  C_{15}\lambda
\end{eqnarray*}
(this calculation also applies to $k=1$).
Thus,
$$|\vphi_k|+\sum_{j} |\vphi_{s_j}| \leq (2C_{14}+C_{15})\,\lambda.$$
If we choose $B=2C_{14}+C_{15}$, \rf{cc5} follows.

Now it is easy to check that \rf{cc6} also holds. Indeed we have
\begin{eqnarray*}
\left(\int_{R_i} |\vphi_i|^p \,d\mu\right)^{1/p} \,\mu(R_i)^{1/p'} & = &
|\alpha_i|\,\mu(A_i) ^{1/p}\,\mu(R_i)^{1/p'} \\
& \leq & C\,|\alpha_i|\,\mu(A_i) \\
& = & C\,\left|\int_{Q_i}f\,w_i\,d\mu\right| \\
& \leq & C\, \left(\int_{Q_i}|f|^p\,d\mu\right)^{1/p}\,\mu(Q_i)^{1/p'},
\end{eqnarray*}
and by \rf{cc1},
\begin{equation} \label{cz47}
\left(\int_{Q_i}|f|^p\,d\mu\right)^{1/p}\,\mu(2Q_i)^{1/p'}
\leq
\frac{C}{\lambda^{p-1}} \int_{Q_i}|f|^p\,d\mu.
\end{equation}
Thus we get \rf{cc6}.

\vspace{3mm}
Suppose now that the collection of cubes $\{Q_i\}_i$ is not finite.
For each fixed $N$ we consider the family of cubes $\{Q_i\}_{1\leq i \leq N}$.
Then, as above, we construct functions $\vphi_1^N,\ldots,\vphi_N^N$ with
$\supp(\vphi_i^N)\subset R_i$ satisfying
$$\int \vphi_i^N \,d\mu = \int_{Q_i} f\,w_i\,d\mu,$$
\begin{equation} \label{cz35}
\sum_{i=1}^N |\vphi_i^N| \leq B\,\lambda
\end{equation}
and, if $1<p<\infty$,
\begin{equation} \label{cz36}
\left(\int_{R_i} |\vphi_i^N|^p \,d\mu\right)^{1/p}\mu(R_i)^{1/p'} \leq
\frac{C}{\lambda^{p-1}} \int_{Q_i}|f|^p\,d\mu.
\end{equation}
By \rf{cz35} and \rf{cz36} there is a subsequence
$\{\vphi_1^k\}_{k\in I_1}$ which is convergent in the weak $\ast$ topology of
$L^\infty(\mu)$ and in the weak $\ast$ topology of $L^p(\mu)$ to some function
$\vphi_1\in L^\infty(\mu)\cap L^p(\mu)$. Now we can consider a subsequence
$\{\vphi_2^k\}_{k\in I_2}$ with $I_2\subset I_1$ which
is convergent also in the weak $\ast$ topologies of
$L^\infty(\mu)$ and $L^p(\mu)$ to some function
$\vphi_2\in L^\infty(\mu)\cap L^p(\mu)$.
In general, for each $j$ we consider a subsequence
$\{\vphi_j^k\}_{k\in I_j}$ with $I_j\subset I_{j-1}$ that converges
in the weak $\ast$ topologies of
$L^\infty(\mu)$ and $L^p(\mu)$ to some function
$\vphi_j\in L^\infty(\mu)\cap L^p(\mu)$.

We have $\supp(\vphi_i)\subset R_i$ and, by the weak $\ast$ convergence in
$L^\infty(\mu)$ and $L^p(\mu)$, the functions $\vphi_i$ also satisfy
\rf{cc4} and \rf{cc6}. To get \rf{cc5}, notice that for each fixed
$m$, by the weak $\ast$ convergence in $L^\infty(\mu)$,
$$\sum_{i=1}^m |\vphi_i| \leq B\lambda,$$
and so \rf{cc5} follows.

\item[c)] For each $i$, we consider the atomic block
$b_i=f\,w_i-\vphi_i$, supported on the cube $R_i$.
Since $K_{Q_i,R_i}\leq C$, by \rf{cz47} and \rf{cc6} we have
$$|b_i|_\hbp \leq \frac{C}{\lambda^{p-1}} \int_{Q_i}|f|^p\,d\mu,$$
which implies \rf{cc7}.
\end{itemize}
\end{proof}


\vspace{5mm}
\begin{proof}[Proof of Theorem \ref{int*}]
For simplicity we assume $\|\mu\|=\infty$.
The proof follows the same lines as the one of \cite[p.43-44]{Journe}.

The functions $f\in L^\infty(\mu)$ having compact support with $\int
f\,d\mu=0$ are dense in $L^{p}(\mu)$, $1<p<\infty$. For such functions we will
show that
\begin{equation}  \label{debpp}
\|M^{\sharp} Tf\|_{L^p(\mu)} \leq C\,\|f\|_{L^p(\mu)}\qquad
1<p<\infty.
\end{equation}
By Theorem \ref{sharp}, this implies
$$\|Tf\|_{L^p(\mu)} \leq C\,\|f\|_{L^p(\mu)}.$$
Notice that if $f\in L^\infty(\mu)$ has compact support and $\int
f\,d\mu=0$, then $f\in\hbm$ and $Tf\in L^1(\mu)$. Thus $N(Tf)\in
L^{1,\infty}(\mu)$, and then $\inf(1,N(Tf))\in L^p(\mu)$. So the
hypotheses of Theorem \ref{sharp} are satisfied.

\vv
Given any function $f\in L^p(\mu)$, $1<p<\infty$, for $\lambda >0$ we take a
family of almost disjoint cubes $\{Q_i\}_i$ as in the previous lemma, and
a collection cubes $\{R_i\}_i$ as in c) in the same lemma. Then we can
write
$$f= b+g = \sum_i \left(\frac{\chi_{Q_i}}{\sum_k \chi_{Q_k}}\,f -
\vphi_i\right) + g.$$
By \rf{cc3} and \rf{cc5}, we have $\|g\|_{L^\infty(\mu)}\leq C\,\lambda$,
and by \rf{cc7},
$$\| b\|_\hbp \leq \frac{C}{\lambda^{p-1}} \|f\|_{L^p(\mu)}^p.$$

Due to the boundedness of $T$ from $L^\infty(\mu)$ into $\rbmo(\mu)$, we have
$$\|M^{\sharp} Tg\|_{L^\infty(\mu)} \leq C_{16}\,\lambda.$$
Therefore,
$$\{M^{\sharp} Tf>(C_{16}+1)\,\lambda \} \subset
\{M^{\sharp} Tb>\lambda \}.$$
Since $M^{\sharp}$ is of weak type $(1,1)$, we have
$$\mu\{M^{\sharp} Tb>\lambda \} \leq
C\,\frac{\|Tb\|_{L^1(\mu)}}{\lambda}.$$
On the other hand, as $T$ is bounded from $\hbm$ into $L^1(\mu)$,
$$\|Tb\|_{L^1(\mu)} \leq C\,\|b\|_\hbm \leq
\frac{C}{\lambda^{p-1}} \|f\|_{L^p(\mu)}^p.$$
Thus
$$\mu\{x\in\R^d:\,M^{\sharp} (Tf)>\lambda \} \leq
C\,\frac{\|f\|_{L^p(\mu)}^p}{\lambda^p}.$$

So the sublinear operator $M^{\sharp} T$ is of weak type $(p,p)$ for all
$p\in (1,\infty)$. By the Marcinkiewitz interpolation theorem
we get that $M^{\sharp} T$ is bounded on $L^p(\mu)$ for all $p\in (1,\infty)$.
In particular, \rf{debpp} holds for a bounded function $f$ with
compact support and $\int f\,d\mu=0$.
\end{proof}


\section{The $T(1)$ theorem for the Cauchy integral}

Before studying the particular case of the Cauchy integral operator,
we will see a result that shows the close relation between CZO's and the
spaces $\rbmo(\mu)$, $\hbm$.

\begin{theorem} \label{equi}
Let $T$ be a CZO and $\rho>1$ some fixed constant. The following
conditions are equivalent:
\begin{itemize}
\item[a)] For any cube $Q$ and any function $a$ supported on $Q$
\begin{equation}  \label{hip111}
\int_Q |T_\ve a|\,d\mu \leq C\,\|a\|_{L^\infty}\,\mu(\rho Q)
\end{equation}
uniformly on $\ve>0$.

\item[b)] $T$ is bounded from $L^\infty(\mu)$ into
$\rbmo(\mu)$.

\item[c)] $T$ is bounded from $\hbm$ into $L^1(\mu)$.
\end{itemize}
\end{theorem}

\begin{proof} We have already seen a)$\Longrightarrow$b) in Theorem \ref{aco1}
and a)$\Longrightarrow$c) in Theorem \ref{aco2}.

Let us prove b)$\Longrightarrow$a). Suppose that $\rho=2$, for example.
Let $a\in L^\infty(\mu)$ be a function
supported on some cube $Q$.
Suppose first $l(Q)\leq \diam(\supp(\mu))/20$ (this is always the case if
$\|\mu\|=\infty$).
We have
\begin{equation} \label{hkk}
\int_Q |T_\ve a - m_\wt{Q} (T_\ve a)|\, d\mu \leq C\,\|a\|_{L^\infty(\mu)}\,\mu(2Q).
\end{equation}
So it is enough to show that
\begin{equation}  \label{hk0}
|m_\wt{Q} (T_\ve a)|\leq C\,\|a\|_{L^\infty(\mu)}.
\end{equation}

Let $x_0\in\supp(\mu)$ be the point (or one of the
points) in $\R^d\setminus (5Q)^\circ$ which is closest to $Q$.
We denote $d_0=\dist(x_0, Q)$.
We assume that $x_0$ is a point
such that some cube with side length $2^{-k}d_0$, $k\geq2$, is
doubling.
Otherwise, we take $y_0$ in $\supp(\mu)\cap B(x_0,l(Q)/100)$ such that
satisfies this condition, and we interchange $x_0$ with $y_0$.

We denote by $R$ a cube concentric with $Q$ with side length
$\max(10d_0,l(\wt{Q}))$. So $K_{\wt{Q},R}\leq C$.
Let $Q_0$ be the biggest doubling cube
centered at $x_0$ with side length $2^{-k}\,d_0$, $k\geq 2$. Then $Q_0\subset
R$, with $K_{Q_0,R}\leq C$, and one can easily check that
\begin{equation}  \label{hk1}
|m_{Q_0} (T_\ve a) - m_\wt{Q} (T_\ve a)|\leq C\,\|T_\ve a\|_{\rbmo(\mu)}\leq
C\,\|a\|_{L^\infty(\mu)}.
\end{equation}
Moreover, $\dist(Q_0,Q)\approx d_0$ and so, for $y\in Q_0$,
$$|T_\ve a(y)| \leq C\,\frac{\mu(Q)}{d_0^n}\,\|a\|_{L^\infty(\mu)}
\leq C\,\|a\|_{L^\infty(\mu)},$$
because $l(Q)<d_0$.
Then we get
$|m_{Q_0} (T_\ve a)|
\leq C\,\|a\|_{L^\infty(\mu)},$
and from \rf{hk1}, we obtain \rf{hk0}.

Suppose now that $l(Q)> \diam(\supp(\mu))/20$. Since $Q$
is centered at some point of $\supp(\mu)$, we may assume that $l(Q)\leq
4\,\diam(\supp(\mu))$. Then $Q\cap \supp(\mu)$ can be covered by a finite
number of cubes $Q_j$ centered at points of $\supp(\mu)$ with side length
$l(Q)/200$. It is quite easy to check that the number of cubes $Q_j$ is
bounded above by some fixed constant $N$ depending only on $d$.
We set
$$a_j = \frac{\chi_{Q_j}}{\sum_k \chi_{Q_k}}\,a.$$
Since a) holds for the cubes $2Q_j$ (which support
the functions $a_j$), we have
\begin{eqnarray*}
\int_Q |T_\ve a|\,d\mu & \leq & \sum_j
\int_{Q\setminus 2Q_j} |T_\ve a_j|\,d\mu
+ \sum_j \int_{2Q_j} |T_\ve a_j|\,d\mu \\
& \leq & \sum_j C\,\|a_j\|_{L^\infty(\mu)} \,\mu(Q) +
\sum_j C\,\|a_j\|_{L^\infty(\mu)} \,\mu(4Q_j) \\
& \leq & C\,N\,\|a\|_{L^\infty(\mu)}\, \mu(2Q).
\end{eqnarray*}

\vspace{2mm}
Now we are going to prove c)$\Longrightarrow$a).
Let $a\in L^\infty(\mu)$ be supported on a cube $Q$.
Assume $\rho=2$ and suppose first
$l(Q)\leq \diam(\supp(\mu))/20$.
We consider the same construction as the one for b)$\Longrightarrow$a).
The cubes $Q$, $Q_0$ and
$R$ are taken as above, and they satisfy $Q,Q_0\subset R$, $K_{Q,R}\leq C$,
$K_{Q_0,R}\leq C$ and $\dist(Q_0,Q) \geq l(Q)$. Recall also that $Q_0$ is
doubling.

We take the atomic block (supported on $R$)
$$b= a + c_{Q_0}\,\chi_{Q_0},$$
where $c_{Q_0}$ is a constant such that $\int b\,d\mu=0$. For $y\in Q$ we have
\begin{equation*}
\begin{split}
|T_\ve(c_{Q_0}\,\chi_{Q_0})(y)| & \leq
C\,\frac{|c_{Q_0}|\,\mu(Q_0)}{\dist(Q,Q_0)^n}  \leq
C\,\frac{\|a\|_{L^1(\mu)}}{\dist(Q,Q_0)^n}\\
& \leq C\,\frac{\,\mu(Q)}{l(Q)^n}  \|a\|_{L^\infty(\mu)}
\leq C\,\|a\|_{L^\infty(\mu)}.
\end{split}
\end{equation*}
Then we have
\begin{eqnarray*}
\int_Q |T_\ve a|\,d\mu & \leq &
\int_Q |T_\ve b|\,d\mu + C\,\|a\|_{L^\infty(\mu)}\,\mu(Q) \\
& \leq & C\,\|b\|_\hbm + C\,\|a\|_{L^\infty(\mu)}\,\mu(Q)  \\
& \leq & C\,K_{Q,R} \,\|a\|_{L^\infty(\mu)}\,\mu(2Q)
+ C\,K_{Q_0,R}\, |c_{Q_0}|\,\mu(2Q_0) \\ \mbox{}
+ C\,\|a\|_{L^\infty(\mu)}\,\mu(Q).
\end{eqnarray*}
Since $Q_0$ is doubling, we have
$$|c_{Q_0}|\,\mu(2Q_0)\leq C\,\|a\|_{L^1(\mu)}
\leq C\,\|a\|_{L^\infty(\mu)}\,\mu(Q).$$
Therefore,
$$\int_Q |T_\ve a|\,d\mu \leq C\,\|a\|_{L^\infty(\mu)}\,\mu(2Q).$$

If $l(Q)> \diam(\supp(\mu))/20$, operating as in the implication
b)$\Longrightarrow$a), we get that a) also holds.
\end{proof}

Now we are going to deal with the $T(1)$ theorem for Cauchy integral operator.
So we take $d=2$ and $n=1$.
Using the relationship of the Cauchy kernel with the curvature of measures,
it is not difficult to get the following result operating
as Melnikov and Verdera \cite{MV}:

\begin{lemma}
Let $\mu$ be some measure on $\C$ satisfying the growth condition
\rf{creix}. If $\|\CE \chi_Q\|_{L^2(\mu_{\mid Q})} \leq C\,\mu(2Q)^{1/2}$
(uniformly on $\ve>0$), then
for any bounded function $a$ with $\supp(a)\subset Q$,
$$\int_Q |\CE a|^2 \,d\mu \leq C\,\|a\|_{L^\infty(\mu)}^2\,\mu(2Q)$$
uniformly on $\ve>0$.
\end{lemma}

We omit the details of the proof (see \cite{MV}). This follows from the
formula, for $a\in L^\infty(\mu)$ with $\supp(a)\subset Q$,
\begin{eqnarray*}
\lefteqn{2\int_Q |\CE a|^2\,d\mu + 4 \real \int_Q a\,\CE a \cdot \overline{\CE
\chi_Q}\, d\mu} && \\
& = &{\uuunt\!\!}_{S_\ve} c(x,y,z)^2\,a(y)\,a(z)\,d\mu(x)\, d\mu(y)\,d\mu(z) +
O(\|a\|_{L^\infty(\mu)}^2\mu(2Q)),
\end{eqnarray*}
where we have denoted
$$S_\ve =\{(x,y,z)\in Q^3:\,|x-y|>\ve,\,|y-z|>\ve,\, |z-x|>\ve\},$$
and $c(x,y,z)$ is the Menger curvature of the triple $(x,y,z)$ (i.e.
the inverse of the radius of the circumference passing through $x,y,z$).

Using the preceeding lemma and the interpolation theorem between the pairs
$(\hbm,L^1(\mu))$ and $(L^\infty(\mu),\rbmo(\mu))$, we get the following version
of the $T1$ theorem for the Cauchy transform for non doubling measures.

\begin{theorem}
\label{T1}
The Cauchy integral operator is bounded on $L^2(\mu)$ if and only if
\begin{equation}  \label{hipt1}
\int_Q |\CE \chi_Q|^2 \,d\mu
\leq C\,\mu(2Q)\quad \mbox{
for any square $Q\subset \C$}
\end{equation}
uniformly on $\ve>0$.
\end{theorem}

This result is already known. The first proofs were obtained independently in
\cite{NTV1} and \cite{Tolsa1}. Another was given later in \cite{Verdera}.
The proof of the present paper follows the lines of the proof
found by Melnikov and Verdera for
the $L^2$ boundedness of the Cauchy integral in the (doubling) case
where $\mu$ is the arc length on a Lipschitz graph \cite{MV}.

Let us remark that in the previous known proofs of the $T(1)$ theorem for the
Cauchy integral, instead of the hypothesis \rf{hipt1}, the assumption was
$$\int_Q |\CE \chi_Q|^2 \,d\mu \leq C\,\mu(Q)
\quad \mbox{for any square $Q\subset \C$.}$$
This is a little stronger than \rf{hipt1}. However, the arguments given
in \cite{NTV1}, \cite{Tolsa1} and \cite{Verdera} can be modified easily
to yield the same result as the one stated in Theorem \ref{T1}.

\vv
Using the relationship between the spaces $\bmo_\rho(\mu)$ and $\rbmo(\mu)$
we obtain another version of the $T(1)$
theorem, which is closer to the classical way of stating the $T(1)$ theorem:

\begin{theorem}  
\label{T1'}
The Cauchy integral operator is bounded on $L^2(\mu)$ if and only if
$\CE(1)\in \bmo_\rho(\mu)$ (uniformly on $\ve>0$), for some $\rho>1$.
\end{theorem}

\begin{proof}
Suppose that $\CE1\in \bmo_\rho(\mu)$. Let us see that this implies
$\CE1\in \rbmo(\mu)$. The estimates are similar to the ones that we used
to show that CZO's bounded on $L^2(\mu)$ are also bounded from $L^\infty(\mu)$
into $\rbmo(\mu)$. Assume, for example $\rho=2$.
We have to show that if $Q\subset R$, then
$$|m_Q(\CE 1) - m_R(\CE 1)|\leq C\,K_{Q,R}\,
\left(\frac{\mu(2Q)}{\mu(Q)} +
\frac{\mu(2R)}{\mu(R)}\right).$$
We denote $Q_R = 2^{N_{Q,R}+1}Q$. Then we
write
\begin{multline*}
|m_Q(\CE 1) - m_R(\CE 1)| \\
\begin{split}
\leq & \,\,|m_Q(\CE \chi_Q)| 
+ |m_Q(\CE\chi_{2Q\setminus Q})| + |m_{Q}(\CE\chi_{Q_R\setminus 2Q})| \\
& + |m_{Q}(\CE\chi_{\C\setminus Q_R}) - m_{R}(\CE\chi_{\C\setminus
Q_R})|\\
& +  |m_R(\CE \chi_R)| + |m_R(\CE\chi_{Q_R\setminus 2R})|
+ |m_R(\CE\chi_{Q_R\cap 2R\setminus R})| \\
= &\,\, M_1 + M_2 + M_3 + M_4 + M_5 + M_6 + M_7.
\end{split}
\end{multline*}
Since $\CE$ is antisymmetric, we have $M_1=M_5=0$. On the other hand, since
the Cauchy transform is bounded from $L^2(\mu_{\mid \C\setminus Q})$ into
$L^2(\mu_{\mid Q})$, we also have 
\begin{eqnarray*}
M_2 & = & |m_Q(\CE\chi_{2Q\setminus Q})| \leq \left(\frac{1}{\mu(Q)}
\int_Q |\CE\chi_{2Q\setminus Q})|^2\, d\mu \right)^{1/2} \\
& \leq &
C\, \left(\frac{\mu(2Q)}{\mu(Q)}\right)^{1/2}
\leq C\,\frac{\mu(2Q)}{\mu(Q)}.
\end{eqnarray*}
By the same argument, we get
$$M_7 = |m_R(\CE\chi_{Q_R\cap 2R\setminus R})| \leq 
C\,\frac{\mu(Q_R\cap 2R)}{\mu(R)} \leq C\,\frac{\mu(2R)}{\mu(R)}.$$
Also, it is easily seen that $|\CE\chi_{Q_R\setminus 2Q}(x)|\leq C\,K_{Q,R}$
for $x\in Q$, and so
$$M_3 = |m_{Q}(\CE\chi_{Q_R\setminus 2Q})| \leq C\,K_{Q,R}.$$
On the other hand, if $x\in Q$ and $y\in R$, we have
$$|\CE\chi_{\C\setminus Q_R}(x) -
\CE \chi_{\C\setminus Q_R}(y)| \leq C,$$
and so $M_4\leq C$. Finally, since $l(Q_R)\approx l(R)$,
$|\CE\chi_{Q_R\setminus 2R}(x)|\leq C$ for $x\in R$, and thus $M_6\leq C$.

Therefore,
\begin{eqnarray*}
|m_Q(\CE 1) - m_R(\CE 1)| & \leq &
C\,K_{Q,R} + C\, \left(\frac{\mu(2Q)}{\mu(Q)} +
\frac{\mu(2R)}{\mu(R)}\right) \\
& \leq &
C\,K_{Q,R} \left(\frac{\mu(2Q)}{\mu(Q)} +
\frac{\mu(2R)}{\mu(R)}\right).
\end{eqnarray*}
So $\CE 1\in \rbmo(\mu)$, and thus we also have $\CE 1\in \bmo_\rho^2(\mu)$, for
any $\rho>1$. Now, some standard calculations show that the condition \rf{hipt1}
of Theorem \ref{T1} is satisfied:
\begin{eqnarray*}
\left(\int_Q |\CE \chi_Q|^2 \,d\mu\right)^{1/2} &=&
\left(\int_Q |\CE \chi_Q - m_Q (\CE \chi_Q)|^2\,d\mu\right)^{1/2} \\
& \leq & \left(\int_Q |\CE 1 - m_Q (\CE 1)|^2 \,d\mu\right)^{1/2}\\
&&\mbox{} +\left(\int_Q |\CE \chi_{\C\setminus 2Q} - 
 m_Q (\CE \chi_{\C\setminus2Q})|^2\, d\mu\right)^{1/2}\\
&&\mbox{} +\left(\int_Q |\CE \chi_{2Q\setminus Q} - 
 m_Q (\CE \chi_{2Q\setminus Q})|^2\, d\mu\right)^{1/2}.
\end{eqnarray*}
Since $\CE 1\in \bmo_2^2(\mu)$, we have
$\int_Q |\CE 1 - m_Q (\CE 1)|^2 \,d\mu \leq C\,\mu(2Q)$.
Also, as usual, we have
$\int_Q |\CE \chi_{\C\setminus 2Q} - 
 m_Q (\CE \chi_{\C\setminus2Q})|^2\, d\mu \leq C\,\,\mu(Q).$
Finally, the last integral can be estimated using the boundedness
of the Cauchy transform from $L^2(\mu_{\mid \C\setminus Q})$ into
$L^2(\mu_{\mid Q})$. Thus \rf{hipt1} holds.  
\end{proof}

Let us remark that, until now, the $T1$ theorem for the Cauchy integral was
known under the assumption $\CE1\in \bmo^2_\rho(\mu)$, but
not under the weaker assumption $\CE1\in \bmo_\rho(\mu)$. 

Also, for general CZO's, 
the assumption $T_\ve1,T^*_\ve1\in \bmo^2_\rho(\mu)$ in the $T1$ theorem 
for non doubling measures of Nazarov, Treil and Volberg
can be substituted by the weaker one $T_\ve1,T^*_\ve1\in \bmo_\rho(\mu)$.
This is due to the fact that if $T_\ve$ is weakly bounded in the sense of
\cite{NTV*} and $T_\ve1,T^*_\ve1\in\bmo_\rho(\mu)$ for some $\rho>1$, then arguing as in
the proof of Theorem \ref{T1'} it follows that $T_\ve1,T^*_\ve1\in\rbmo(\mu)$, and so
$T_\ve1,T^*_\ve1\in \bmo^2_\rho(\mu)$, for any $\rho>1$.


\section{Commutators}

In this section we will prove that if $b\in\rbmo(\mu)$ and $T$ is a CZO
bounded on $L^2(\mu)$, then the commutator $[b,T]$ defined by
$$[b,T](f) =b\,T(f) - T(bf)$$
is bounded on $L^p(\mu)$, $1<p<\infty$.
In this formula, $T$ stands for a weak limit as
$\ve\to0$ of some subsequence of the uniformly bounded operators $T_\ve$.

The $L^p(\mu)$ boundedness of the commutator $[b,T]$
is a result due to Coifman,
Rochberg and Weiss \cite{CRW} in the classical case where $\mu$ is the
Lebesgue measure on $\R^d$.
Their proof, with some minor changes, works also
for doubling measures.
On the other hand, for $\mu$ being the Lebesgue measure, they showed that if
$R_i$, $i=1,\cdots,d$, are the Riesz transforms on $\R^d$, then the
$L^p(\mu)$ boundedness of the
commutators $[b,R_i]$, $i=1,\cdots,d$, for some $p\in(1,\infty)$
implies $b\in BMO(\mu)$.

When $\mu$ is a non doubling measure and $b$ satisfies \rf{dbmo}, i.e. it
belongs to the classical space $\bmo(\mu)$, then it has been shown by
Orobitg and P\'erez \cite{OP} that the commutator $[b,T]$ is bounded on
$L^p(\mu)$, $1<p<\infty$.

Let us state now the result that we will obtain in this section in detail.

\begin{theorem} \label{commut}
If $T$ is a CZO bounded on $L^2(\mu)$ and $b\in \rbmo(\mu)$, then the
commutator $[b,T]$ is bounded on $L^p(\mu)$.
\end{theorem}

Our proof will be based on the use of
the sharp maximal operator, as the one of Janson and Str\"omberg \cite{Janson}
for the doubling case. However, the result can be obtained also by means
of a good $\lambda$ inequality, as in \cite{CRW}.

To prove Theorem \ref{commut} we will need a couple of lemmas dealing with the
coefficients $K_{Q,R}$.

\begin{lemma} \label{pok0}
There exists some constant $P$ (big enough) depending on $C_0$ and $n$ such
that if $Q_1\subset Q_2\subset \cdots \subset Q_m$ are concentric cubes with
$K_{Q_i,Q_{i+1}} >P$ for $i=1,\ldots,m-1$, then
\begin{equation}  \label{uyy0}
\sum_{i=1}^{m-1} K_{Q_i,Q_{i+1}} \leq C_{17}\, K_{Q_1,Q_m},
\end{equation}
where $C_{17}$  depends only on $C_0$ and $n$.
\end{lemma}

\begin{proof} Let $Q_i'$ be a cube concentric with $Q_i$ such that
$l(Q_i)\leq l(Q_i') < 2l(Q_i)$, with $l(Q_i')=2^k l(Q_1)$ for some $k\geq0$.
Then
$$C_{18}^{-1}\,K_{Q_i,Q_{i+1}}\leq K_{Q_i',Q_{i+1}'}\leq
C_{18}\,K_{Q_i,Q_{i+1}},$$
for all $i$, with $C_{18}$ depending on $C_0$ and $n$.

Observe also that if we take $P$ so that $C_{18}^{-1}P\geq 2$, then
$K_{Q_i',Q_{i+1}'}>2$ and so
$$K_{Q_i',Q_{i+1}'} \leq 2 \sum_{k=1}^{N_{Q_i',Q_{i+1}'}}
\frac{\mu(2^k Q_i')}{l(2^kQ_i')^n}.$$
Therefore,
\begin{equation} \label{uyy1}
\sum_i K_{Q_i',Q_{i+1}'} \leq 2 \sum_i \sum_{k=1}^{N_{Q_i',Q_{i+1}'}}
\frac{\mu(2^k Q_i')}{l(2^kQ_i')^n}.
\end{equation}
On the other hand, if $P$ is big enough, then $Q_i'\neq Q_{i+1}'$. Indeed,
$$C_0\,N_{Q_i,Q_{i+1}} \geq \sum_{k=1}^{N_{Q_i,Q_{i+1}}} \frac{\mu(2^k Q_i)}{
l(2^k Q_i)^n} \geq P-1,$$
and so $N_{Q_i,Q_{i+1}} \geq (P-1)/C_0>2$, assuming $P$ big enough. This
implies $l(Q_{i+1})>2l(Q_i)$, and then, by construction, $Q_i'\neq Q_{i+1}'$.

As a consequence, on the right hand side of \rf{uyy1}, there is no overlapping
in the terms $\frac{\mu(2^k Q_i')}{l(2^k Q_i')^n}$, and then
$$\sum_i K_{Q_i',Q_{i+1}'} \leq 2 K_{Q_1,Q_m'} \leq 2C_{18}\, K_{Q_1,Q_m},$$
and \rf{uyy0} follows.
\end{proof}

\begin{lemma} \label{pok}
There exists some constant $P_0$ (big enough) depending on $C_0$, $n$ and
$\beta_d$ such that if $x\in \R^d$ is some fixed point and $\{f_Q\}_{Q\ni x}$ 
is a collection of numbers
such that $|f_Q-f_R|\leq C_x$ for all doubling cubes $Q\subset R$ with $x\in
Q$ such that $K_{Q,R} \leq P_0$, then
$$|f_Q-f_R|\leq C\,K_{Q,R}\,C_x \quad \mbox{for all doubling
cubes $Q\subset R$ with $x\in Q$},$$
where $C$ depends on $C_0$, $n$, $P_0$ and $\beta_d$.
\end{lemma}

\begin{proof} Let $Q\subset R$ be two doubling cubes in $\R^d$, with $x\in Q =:
Q_0$. Let $Q_1$ be the first cube of the form $2^k Q$, $k\geq0$, such that
$K_{Q,Q_1}>P$. Since $K_{Q,2^{-1}Q_1}\leq P$, we have $K_{Q,Q_1}\leq P+C_0$.
Therefore, for the doubling cube $\wt{Q}_1$, we have $K_{Q,\wt{Q}_1}\leq
C_{19}$, with $C_{19}$ depending on $P$, $n$, $\beta_d$ and $C_0$.

In general, given $\wt{Q}_i$, we denote by $Q_{i+1}$ the first cube
of the form $2^k \wt{Q}_i$, $k\geq0$, such that $K_{\wt{Q}_i,Q_{i+1}}>P$,
and we consider the cube $\wt{Q}_{i+1}$. Then, we have
$K_{\wt{Q}_i,\wt{Q}_{i+1}}\leq C_{19}$, and also
$K_{\wt{Q}_i,\wt{Q}_{i+1}}>K_{\wt{Q}_i,Q_{i+1}} > P$.

Then we obtain
\begin{equation} \label{uyy3}
|f_Q-f_R|\leq \sum_{i=1}^N |f_{\wt{Q}_{i-1}} - f_{\wt{Q}_{i}}| +
|f_{\wt{Q}_N} - f_R|,
\end{equation}
where $\wt{Q}_N$ is the first cube of the sequence $\{\wt{Q}_i\}_i$ such that
$\wt{Q}_{N+1}\supset R$. Since
$K_{\wt{Q}_N,\wt{Q}_{N+1}}\leq C_{19}$, we also have
$K_{\wt{Q}_N,R}\leq C_{19}$.
By \rf{uyy3} and Lemma \ref{pok0}, if we set $P_0=C_{19}$, we get
\begin{eqnarray*}
|f_Q-f_R| & \leq & \sum_{i=1}^N K_{\wt{Q}_i,\wt{Q}_{i+1}} \,C_x +
K_{\wt{Q}_N,R}\,C_x \\
& \leq & C\,K_{Q,\wt{Q}_{N}} \,C_x + K_{\wt{Q}_N,R}\,C_x
\leq C\,K_{Q,R} \,C_x.
\end{eqnarray*}
\end{proof}

\vspace{5mm}
\brem
By the preceeding lemma, to see if some function $f$ belongs to $\rbmo(\mu)$,
the regularity condition \rf{def2} only needs to be checked for doubling
cubes $Q\subset R$ such that $K_{Q,R} \leq P_0$. In a similar way, it can be
proved that if the regularity condition \rf{def2'} holds for any pair of 
cubes $Q\subset R$ with $K_{Q,R}$ not too large, then it holds for any pair of 
cubes $Q\subset R$.

On the other hand, one can introduce an operator ${\wh{M}^\sharp}$ defined as
$M^\sharp$, but with the second supremum in the definition \rf{defsharp}
taken only over doubling cubes $Q\subset R$ such that $x\in Q$ and
$K_{Q,R}\leq P_0$. Then, by the preceeding lemma it follows that
${\wh{M}^\sharp}(f)\approx M^\sharp(f)$.
\erem


\vspace{5mm}
\begin{proof}[Proof of Theorem \ref{commut}]
For all $p\in(1,\infty)$, we will show the
pointwise inequality
\begin{equation} \label{pwi}
M^\sharp([b,T]f)(x) \leq C_p\,\|b\|_*\,(M_{p,(9/8)} f(x) +
M_{p,(3/2)} T f(x) + T_* f(x)),
\end{equation}
where, for $\eta>1$, $M_{p,(\eta)}$ is the non centered maximal operator
$$M_{p,(\eta)} f(x) = \sup_{Q\ni x}  \left(\frac{1}{\mu(\eta Q)}
\int_Q |f|^p\,d\mu\right)^{1/p},$$
and $T_*$ is defined as
$$T_*f(x) = \sup_{\ve>0} |T_{\ve}f(x)|.$$
The operator $M_{p,(\eta)}$ is bounded on $L^r(\mu)$ for $r>p$, and $T_*$
is bounded on $L^r(\mu)$ for $1<r<\infty$ because $T$ is bounded on $L^2(\mu)$
(see \cite{NTV2}). Then
the pointwise inequality \rf{pwi} for $1<p<\infty$ implies
the $L^p(\mu)$ boundedness of $M^\sharp([b,T])$ for $1<p<\infty$.
If $b$ is a bounded
function we can apply Theorem \ref{sharp} because, by the $L^p(\mu)$
boundedness of $T$, it follows that $[b,T]\in L^p(\mu)$. On the other hand,
by Lemma \ref{lemaq} it is easily seen that we can assume that
$b$ is a bounded function. So the inequality \rf{pwi} implies that $[b,T]$
is bounded on $L^p(\mu)$, $1<p<\infty$.

Let $\{b_Q\}_Q$ a family of numbers satisfying
$$\int_Q |b-b_Q|\,d\mu\leq 2\mu(2Q)\,\|b\|_{**}$$
for any cube $Q$, and $$|b_Q-b_R|\leq 2K_{Q,R}\,\|b\|_{**}$$ for all cubes
$Q\subset R$. For any cube $Q$, we denote
$$h_Q := m_Q(T ((b-b_Q)\,f\chi_{\R^d\setminus \tfrac{4}{3}Q})$$
We will show that
\begin{equation} \label{pw1}
\frac{1}{\mu(\ts \frac{3}{2} Q)} \int_Q |[b,T]f - h_Q| \,d\mu
\leq C\,\|b\|_*\,(M_{p,(9/8)} f(x) +
M_{p,(3/2)} T f(x))
\end{equation}
for all $x$ and $Q$ with $x\in Q$, and
\begin{equation} \label{pw2}
|h_Q- h_R| \leq C\,\|b\|_*\,(M_{p,(9/8)} f(x) + T_* f(x))\,K_{Q,R}^2
\end{equation}
for {\em all} cubes $Q\subset R$ with $x\in Q$.
In the final part of the proof we will see that from the preceeding two
inequalities one easily gets \rf{pwi}.

To get \rf{pw1} for some fixed cube $ Q$ and $x$ with $x\in Q$,
we write $[b,T]f$ in the following way:
\begin{eqnarray} \label{pw3}
[b,T]f & = & (b-b_Q)\,T f - T((b-b_Q)\,f) \nonumber \\
& = & (b-b_Q)\,T f - T((b-b_Q)\,f_1) - T((b-b_Q)\,f_2),
\end{eqnarray}
where $f_1= f\,\chi_{\frac{4}{3}Q}$ and $f_2=f-f_1$.
Let us estimate the term $(b-b_Q)\,T f$:
\begin{eqnarray} \label{pw4}
\frac{1}{\mu(\ts \frac{3}{2}Q)} \int_Q |(b-b_Q)\,T f|\, d\mu & \leq &
\left(\frac{1}{\mu(\ts \frac{3}{2}Q)}
\int_Q |(b-b_Q)|^{p'}\, d\mu\right)^{1/p'} \nonumber \\
& & \times
\left(\frac{1}{\mu(\ts \frac{3}{2}Q)}
\int_Q |T f|^p\, d\mu\right)^{1/p}\nonumber \\
&&\nonumber \\
& \leq & C\,\|b\|_*\,M_{p,(3/2)}Tf(x).
\end{eqnarray}
Now we are going to estimate the second term on the right hand side of
\rf{pw3}. We take $s=\sqrt{p}$. Then we have
\begin{eqnarray*}
\left[\frac{1}{\mu(\ts \frac{3}{2}Q)} \int_{\frac{4}{3}Q}
 |(b-b_Q)\,f_1|^s\, d\mu \right]^{1/s}
& \leq & \left(\frac{1}{\mu(\ts \frac{3}{2}Q)}
\int_{\frac{4}{3}Q} |b-b_Q|^{ss'}\, d\mu\right)^{1/ss'} \\
&& \times
\left(\frac{1}{\mu(\ts \frac{3}{2}Q)}
\int_{\frac{4}{3}Q} |f|^p\, d\mu\right)^{1/p}  \\ &&\\
& \leq & C\,\|b\|_*\,M_{p,(9/8)}f(x).
\end{eqnarray*}
Notice that we have used that
$\int_{\frac{4}{3}Q} |b-b_Q|^{ss'}\,d\mu \leq
C\,\|b\|_*^{ss'}\,\mu(\frac{3}{2} Q),$
which holds because $|b_Q- b_{\frac{4}{3}Q}| \leq C\,\|b\|_*$.
Then we get
\begin{eqnarray}   \label{pw5}
\frac{1}{\mu(\ts \frac{3}{2}Q)} \int_Q |T((b-b_Q)\,f_1)|\,d\mu & \leq &
\frac{\mu(Q)^{1-1/s}}{\mu(\ts \frac{3}{2}Q)}\,
\|T((b-b_Q)\,f_1)\|_{L^s(\mu)} \nonumber \\
& \leq &
C\, \frac{\mu(Q)^{1-1/s}}{\mu(\ts \frac{3}{2}Q)}\, \|(b-b_Q)\,f_1\|_{L^s(\mu)}
\nonumber \\
& \leq &
C\,\frac{\mu(Q)^{1-1/s}}{\mu(\ts \frac{3}{2}Q)^{1-1/s}}\,
\|b\|_*\,M_{p,(9/8)}f(x) \nonumber \\ & \leq & C\, \|b\|_*\,M_{p,(9/8)}f(x).
\end{eqnarray}

By \rf{pw3}, \rf{pw4} and \rf{pw5}, to prove \rf{pw1}
we only have to estimate the difference $|T((b-b_Q)\,f_2) - h_Q|$. For 
$x,y\in Q$ we have
\begin{multline} \label{nono}
|T((b-b_Q)\,f_2)(x) -
T((b-b_Q)\,f_2)(y)| \\
\begin{split}
& \leq
C \int_{\R^d\setminus \frac{4}{3}Q} \frac{|y-x|^\delta}{|z-x|^{n+\delta}}\,
|b(z)-b_Q| \,|f(z)|\,d\mu(z) \\
& \leq  C \sum_{k=1}^\infty
\int_{2^k\frac{4}{3}Q\setminus 2^{k-1}\frac{4}{3}Q}
\frac{l(Q)^\delta}{|z-x|^{n+\delta}}\,
(|b(z)-b_{2^k\frac{4}{3}Q}| + |b_Q-b_{2^k\frac{4}{3}Q}|) \,|f(z)|\,d\mu(z) \\
& \leq  C \sum_{k=1}^\infty 2^{-k\delta} \frac{1}{l(2^kQ)^n}
\int_{2^k\frac{4}{3}Q} |b(z)- b_{2^k\frac{4}{3}Q}|\,|f(z)| \,d\mu(z) \\
& +
C \sum_{k=1}^\infty k2^{-k\delta}\,\|b\|_*\, \frac{1}{l(2^kQ)^n}
\int_{2^k\frac{4}{3}Q} |f(z)|\, d\mu(z) \\
& \leq  C \sum_{k=1}^\infty 2^{-k\delta}\,\|b\|_*\,M_{p,(9/8)}f(x) +
C \sum_{k=1}^\infty k2^{-k\delta}\,\|b\|_*\,M_{(9/8)}f(x) \\
& \leq  C\, \|b\|_*\,M_{p,(9/8)}f(x).
\end{split}
\end{multline}
Taking the mean over $y\in Q$, we get
\begin{eqnarray*}
|T((b-b_Q)\,f_2)(x) - h_Q| & = &|T((b-b_Q)\,f_2)(x) - m_Q(T((b-b_Q)\,f_2))| \\
& \leq & C\, \|b\|_*\,M_{p,(9/8)}f(x),
\end{eqnarray*}
and so \rf{pw1} holds.

\vspace{4mm}
Now we have to check the regularity condition \rf{pw2} for the numbers
$\{h_Q\}_Q$. Consider two cubes $Q\subset R$ with $x\in Q$. We denote
$N=N_{Q,R}+1$. We write the difference $|h_Q-h_R|$
in the following way:
\begin{multline*}
|m_Q(T((b-b_Q)\,f\chi_{\tfrac{4}{3}Q})) - 
m_R(T((b-b_R)\,f\chi_{\tfrac{4}{3}R}))| \\
\begin{split}
\leq &\,\, |m_Q(T((b-b_Q)\,f\chi_{2Q\setminus\tfrac{4}{3}Q}))| \\
& +|m_Q(T((b_Q-b_R)\,f\,\chi_{\R^d\setminus2Q}))| \\
&+ |m_Q(T((b-b_R)\,f\,\chi_{2^N Q\setminus2Q}))| \\
&+ |m_Q(T((b-b_R)\,f\,\chi_{\R^d\setminus2^N Q})) -
m_R(T((b-b_R)\,f\,\chi_{\R^d\setminus2^NQ}))| \\
&+ |m_R(T((b-b_R)\,f\,\chi_{2^N Q\setminus2R}))| \\
= & \,\,M_1+ M_2 + M_3 + M_4+M_5.
\end{split}
\end{multline*}
Let us estimate $M_1$. For $y\in Q$ we have
\begin{eqnarray*}
|T((b-b_Q)\,f\,\chi_{2Q\setminus\tfrac{4}{3}Q})(y)| & \leq &
\frac{C}{l(Q)^n}\int_{2Q}|b-b_Q|\, |f|\,d\mu  \\
& \leq & \frac{C}{l(Q)^n} \left(\int_{2Q}|b-b_Q|^{p'}\,d\mu
\right)^{1/p'}\,
\left(\int_{2Q}|f|^{p}\,d\mu\right)^{1/p} \\
& \leq & \frac{C\,\|b\|_*}{l(Q)^{n/p}} \left(\int_{2Q}|f|^{p}\,d\mu\right)^{1/p}\\
& \leq & C\,\|b\|_*\,M_{p,(9/8)}f(x).
\end{eqnarray*}
So we derive $M_1\leq  C\,\|b\|_*\,M_{p,(9/8)}f(x)$.

Let us consider the term $M_2$. For $x,y\in Q$, it is easily seen that
\begin{eqnarray*}
|T(f\,\chi_{\R^d\setminus2Q})(y)|
& \leq & T_*f(x) + C\,\sup_{{Q_0}\ni x} \frac{1}{l(Q_0)^n} \int_{Q_0} |f|\,d\mu
\\ & \leq & T_*f(x) + C\,M_{p,(9/8)}f(x).
\end{eqnarray*}
Thus
$$M_2 = |(b_R-b_Q)\,T(f\,\chi_{\R^d\setminus2Q})(y)|
\leq C\,K_{Q,R}\,\|b\|_*\,(T_*f(x) + M_{p,(9/8)}f(x)).$$

Let us turn our attention to the term $M_4$. Operating as in \rf{nono}, for any 
$y,z\in R$, we get
$$|T((b-b_R)\,f\,\chi_{\R^d\setminus2^NQ})(y) -
T((b-b_R)\,f\,\chi_{\R^d\setminus2^NQ})(z)|
\leq  C\, \|b\|_*\,M_{p,(9/8)}f(x).$$
Taking the mean over $Q$ for $y$ and over $R$ for $z$, we obtain
$$M_4\leq  C\, \|b\|_*\,M_{p,(9/8)}f(x).$$
The term $M_5$ is easy to estimate too. Some calculations very similar to the
ones for $M_1$ yield $M_5\leq  C\,\|b\|_*\,M_{p,(9/8)}f(x)$.

Finally, we have to deal with $M_3$.
For $y\in Q$, we have
\begin{eqnarray*}
|T((b-b_R)\,f\,\chi_{2^N Q\setminus2Q}))(y)| & \leq &
C\, \sum_{k=1}^{N-1}  \frac{1}{l(2^kQ)^n} \int_{2^{k+1}Q\setminus 2^kQ}
|b-b_R|\,|f|\,d\mu \\
& \leq & C\, \sum_{k=1}^{N-1} \frac{1}{l(2^kQ)^n}
\left(\int_{2^{k+1}Q} |b-b_R|^{p'}\,d\mu \right)^{1/p'} \\
&&\times
\left(\int_{2^{k+1}Q} |f|^p\,d\mu \right)^{1/p}.
\end{eqnarray*}
We have
\begin{eqnarray*}
\left(\int_{2^{k+1}Q} |b-b_R|^{p'}\,d\mu \right)^{1/p'} & \leq &
\left(\int_{2^{k+1}Q} |b-b_{2^{k+1}Q}|^{p'}\,d\mu \right)^{1/p'} \\ && \mbox{}+
\mu(2^{k+1}Q)^{1/p'}\,|b_{2^{k+1}Q}-b_R| \\
&&\\
& \leq & C\,K_{Q,R}\,\|b\|_*\,\mu(2^{k+2}Q)^{1/p'}.
\end{eqnarray*}
Thus
\begin{multline*}
|T((b-b_R)\,f\,\chi_{2^N Q\setminus2Q}))(y)| \\
\begin{split} & \leq
C\,K_{Q,R}\,\|b\|_*
 \sum_{k=1}^{N-1} \frac{\mu(2^{k+2}Q)}{l(2^kQ)^n}
\left(\frac{1}{\mu(2^{k+2}Q)}\int_{2^{k+1}Q} |f|^p\,d\mu \right)^{1/p}\\
& \leq  C\,K_{Q,R}\,\|b\|_*\, M_{(9/8)}f(x)
 \sum_{k=1}^{N-1} \frac{\mu(2^{k+2}Q)}{l(2^kQ)^n} \\
& \leq  C\,K_{Q,R}^2\,\|b\|_*\, M_{p,(9/8)}f(x).
\end{split}
\end{multline*}
Taking the mean over $Q$, we get
$$M_3 \leq C\,K_{Q,R}^2\,\|b\|_*\, M_{p,(9/8)}f(x).$$
So by the estimates on $M_1,M_2,M_3,M_4$ and $M_5$, the regularity condition
\rf{pw2} follows.

\vspace{4mm}
Let us see how from \rf{pw1} and \rf{pw2} one obtains \rf{pwi}.
 From \rf{pw1}, if $Q$ is a {\em doubling} cube and $x\in Q$, we have
\begin{eqnarray}  \label{pw98}
|m_Q([b,T]f) - h_Q| & \leq & \frac{1}{\mu(Q)} \int_Q
|[b,T]f - h_Q| \,d\mu \nonumber \\
& \leq & C\,\|b\|_*\,(M_{p,(9/8)} f(x) + M_{p,(3/2)}
Tf(x)).
\end{eqnarray}
Also, for any cube $Q\ni x$ (non doubling, in general),
$K_{Q,\wt{Q}}\leq C$, and then by \rf{pw1} and \rf{pw2} we get
\begin{multline} \label{pw99}
\frac{1}{\mu(\ts \frac{3}{2} Q)} \int_Q |[b,T]f - m_\wt{Q}([b,T]f)
| \,d\mu \\
\begin{split} & \leq
\frac{1}{\mu(\ts \frac{3}{2} Q)} \int_Q |[b,T]f - h_Q| \,d\mu +
|h_Q-h_\wt{Q}| + |h_\wt{Q}-m_\wt{Q}([b,T]f)|  \\
& \leq  C\,\|b\|_*\,(M_{p,(9/8)} f(x) + M_{p,(3/2)} Tf(x)+ T_* f(x)).
\end{split}
\end{multline}

On the other hand,
for all {\em doubling} cubes $Q\subset R$ with $x\in Q$ such that $K_{Q,R}\leq
P_0$, where
$P_0$ is the constant in Lemma \ref{pok}, by \rf{pw2} we have
$$|h_Q- h_R| \leq C\,\|b\|_*\,(M_{p,(9/8)} f(x) + T_* f(x))\,P_0^2.$$
So by Lemma \ref{pok} we get
$$|h_Q- h_R| \leq C\,\|b\|_*\,(M_{p,(9/8)} f(x) + T_* f(x))\,K_{Q,R}$$
for all doubling cubes $Q\subset R$ with $x\in Q$
and, using \rf{pw98} again, we obtain
\begin{multline*}
|m_Q([b,T]f) - m_R([b,T]f)| \\ \leq
C\,\|b\|_*\,(M_{p,(9/8)} f(x) + M_{p,(3/2)} Tf(x) + T_* f(x))\,K_{Q,R}.
\end{multline*}
 From this estimate and \rf{pw99}, we get \rf{pwi}. \end{proof}


\vspace{6mm}
\noindent {\bf Acknowledgements.}
The author wishes to thank the hospitality of the
Department of Mathematics of the University of Jyv\"{a}skyl\"{a} where part of this research
was done. Thanks are due especially to Pertti Mattila.

The author is also grateful to the referee for carefully reading the paper
and for helpful comments.
\vspace{10mm}



\begin{thebibliography}{MMNO}

\bibitem[CRW]{CRW} R. Coifman, R. Rochberg and G. Weiss. {\em Factorization
theorems for Hardy spaces in several variables.} Ann. of Math. 103 (1976),
611-635.

\bibitem[Da]{David} G. David. {\em Unrectifiable $1$-sets have vanishing
analytic capacity.} Revista Mat. Iberoamericana 14:2 (1998), 369-479.

\bibitem[GR]{GR} J. Garc\'{\i}a-Cuerva and J.L. Rubio de Francia. {\em Weighted norm
inequalities and related topics.} North-Holland Math. Studies 116, 1985.

\bibitem[Gu]{Guzman} M. de Guzm\'an. {\em Differentiation of integrals in
$\R^n$.} Lecture Notes in Math. 481, Springer-Verlag, 1975.

\bibitem[Ja]{Janson} S. Janson. {\em Mean oscillation and commutators of
singular integral operators.} Ark. Mat. 16 (1978), 263-270.

\bibitem[Jo]{Journe} J.-L. Journ\'e. {\em Calder\'on-Zygmund operators,
pseudo-differential operators and the Cauchy integral of Calder\'on.}
Lecture Notes in Math. 994, Springer-Verlag, 1983.

\bibitem[MMNO]{MMNO} J. Mateu, P. Mattila, A. Nicolau and J. Orobitg. {\em
$\bmo$ for non doubling measures.} To appear in Duke Math. J.

\bibitem[MV]{MV} M. S. Melnikov and J. Verdera. {\em A geometric proof
of the $L^2$ boundedness of the Cauchy integral on Lipschitz graphs.} Int.
Math. Res. Not. {\bf 7} (1995), 325-331.

\bibitem[Mo]{Morse} A.P. Morse. {\em Perfect blankets.} Trans. Amer. Math. Soc.
{\bf6} (1947), 418-442.

\bibitem[NTV1]{NTV1}  F. Nazarov, S. Treil and A. Volberg. {\em Cauchy integral
and Calder\'on-Zygmund operators on nonhomogeneous spaces.}
Int. Math. Res. Not. {\bf 15} (1997), 703-726.

\bibitem[NTV2]{NTV2} {F. Nazarov, S. Treil and A. Volberg.}
{\em Weak type estimates and Cotlar inequalities for
Calder\'on-Zygmund operators in nonhomogeneous spaces.}
Int. Math. Res. Not. {\bf 9} (1998), 463-487.

\bibitem[NTV3]{NTV3} F. Nazarov, S. Treil and A. Volberg,
{\em Pulling ourselves up by the hair.} Preprint (1997).

\bibitem[NTV4]{NTV*}
F. Nazarov, S. Treil and A. Volberg. {\em $Tb$-theorem on non-homogeneous
spaces.} Preprint (1999).

\bibitem[OP]{OP} J. Orobitg and C. P\'erez. {\em Some analysis for non doubling
measures.} Preprint (1999).

\bibitem[To1]{Tolsa1} X. Tolsa. {\em $L^2$-boundedness of the Cauchy integral
operator for continuous measures.} Duke Math. J. {\bf 98}:2 (1999), 269-304.

\bibitem[To2]{Tolsa2} X. Tolsa. {\em Cotlar's inequality and existence of
principal values for the Cauchy integral without the doubling condition.}
J. Reine Angew. Math. {\bf 502} (1998), 199-235.

\bibitem[To3]{Tolsa3} X. Tolsa. {\em A $T(1)$ theorem for non doubling
measures with atoms.} To appear in Proc. London Math. Soc.

\bibitem[Ve]{Verdera} J. Verdera. {\em On the $T(1)$ theorem for the
Cauchy integral.} To appear in Ark. Mat.

\bibitem[Yo]{Yosida} K. Yosida. {\em Functional Analysis.}
Springer-Verlag ($6^{th}$ ed.), 1980.

\end{thebibliography}
\end{document}